\documentclass[12pt, requno]{amsart}
\usepackage{amsmath}
\usepackage{amssymb}
\usepackage{epsfig}
\usepackage{graphicx}
\usepackage{color}
\definecolor{shadecolor}{gray}{0.875}
\usepackage{amscd}
\usepackage{comment}
\usepackage{hyperref}
\usepackage{tikz-cd}

\numberwithin{equation}{section}

\input xy
\xyoption{all}

\calclayout
\allowdisplaybreaks[3]

\theoremstyle{plain}
\newtheorem{prop}{Proposition}[section]

\newtheorem{theo}[prop]{Theorem}
\newtheorem{coro}[prop]{Corollary}

\newtheorem{lemm}[prop]{Lemma}

\theoremstyle{definition}
\newtheorem{defi}[prop]{Definition}

\newtheorem{ques}[prop]{Question}

\newtheorem{rema}[prop]{Remark}

\def\bQ{{\mathbb Q}}

\def\bZ{{\mathbb Z}}

\def\Nef{\mathrm{Nef}}
\def\Hilb{\mathrm{Hilb}}

\def\Sec{\mathrm{Sec}}

\makeatother
\makeatletter

\author{Sho Tanimoto}
\address{Graduate School of Mathematics, Nagoya University, Furocho Chikusa-ku, Nagoya, 464-8602, Japan}
\email{sho.tanimoto@math.nagoya-u.ac.jp}

\title[The spaces of rational curves on del Pezzo surfaces]{The spaces of rational curves \\on del Pezzo surfaces via conic bundles}

\begin{document}
%\date{\today}

\maketitle

{
\centering
{\it Dedicated to Yuri Tschinkel on his 60th birthday}\par
}

\begin{abstract}
Using the homological sieve method developed by Das, Lehmann, Tosteson, and the author, we prove Peyre's all height approach to Manin's conjecture for split quintic del Pezzo surfaces defined over $\mathbb F_q(t)$ assuming $q$ is sufficiently large. We also establish lower bounds of correct magnitude for the counting function of rational curves on split low degree del Pezzo surfaces defined over $\mathbb F_q$ assuming $q$ is large.
\end{abstract}

\tableofcontents

\section{Introduction}

One of the central questions in diophantine geometry is the distribution of rational points on algebraic varieties, i.e., the density of the set of rational points in various topologies or the asymptotic formula for the counting function of rational points. Around the 1980s, Yuri Manin and his collaborators have started an ambitious program to formulate the asymptotic formula for the counting function of rational points on a smooth Fano variety defined over a number field. This lead to the formulation of Manin's conjecture which has evolved over the last three decades. (\cite{FMT, BM, Peyre, BT, Peyre03, Peyre17, LST18, LS24})

One can formulate Manin's conjecture over global function fields, and there is an influential heuristic on this conjecture by Victor Batyrev (\cite{Bat88}). An idea is that counting rational points over global function fields amounts to counting the number of $\mathbb F_q$-points on the moduli space of curves/sections, and Batyrev suggested to use the geometry of the moduli spaces to prove Manin's conjecture over such fields. This has been further strengthened by Jordan Ellenberg and Akshay Venkatesh who proposed to use the homological stability of the moduli space and combine it with the Grothendieck--Lefschetz trace formula to prove the conjecture. (\cite{EV05}) In \cite{DLTT25}, Das, Lehmann, Tosteson, and the author developed the homological sieve method, a geometric/topological approach to formalize the inclusion-exclusion principle and used this to prove a version of Manin's conjecture for rational curves on split quartic del Pezzo surfaces defined over $\mathbb F_q$ with $q$ being large. This work was based on Das--Tosteson's bar complex method which has been developed to establish homological stability of these moduli spaces. (\cite{DT24})

In this paper, we discuss implications of the homological sieve method developed in \cite{DLTT25} for other split del Pezzo surfaces defined over finite fields.

\subsection{Main results}

Let $S$ be a smooth split del Pezzo surface defined over $k = \mathbb F_q$.
Let $\mathrm{Nef}_1(S)$ be the nef cone of $S$. 
For each nef class $\alpha \in \mathrm{Nef}_1(S)_{\mathbb Z}$, let $\mathrm{Mor}(\mathbb P^1, S, \alpha)$ be the morphism scheme parametrizing rational curves $f : \mathbb P^1 \to S$ whose class is $\alpha$. Let $\mathcal C \subset \mathrm{Nef}_1(S)$ be a rational polyhedral cone. We define the counting function as
\[
N(\mathbb P^1, S, -K_S, \mathcal C, d) = \sum_{\alpha \in \mathcal C_{\mathbb Z}, -K_S.\alpha \leq d} \#\mathrm{Mor}(\mathbb P^1, S, \alpha)(k).
\]
This number is equal to the number of $\mathbb F_q$-rational curves $f : \mathbb P^1 \to S$ of anticanonical degree $\leq d$, whose class is contained in $\mathcal C$.

\subsubsection*{Quintic del Pezzo surfaces}

Suppose that $S$ is a smooth split quintic del Pezzo surface over $k$, i.e., $(-K_S)^2 = 5$ and $\rho(S) = \rho(S_{\overline{k}})$.
Here is one of our results: 
\begin{theo}
\label{theo:intromainI}
Let $S$ be a smooth split quintic del Pezzo surface over $k = \mathbb F_q$.
Then there exist uniform constants $C \geq 1$ and $C_3 \geq 1$ independent of $S$ and a rational homogeneous continuous piecewise linear function $\ell(\alpha)$ on $\mathrm{Nef}_1(S)$ which is independent of $S$ and which takes positive values on the interior $\mathrm{Nef}_1(S)^\circ$ with the following properties: let $\alpha$ be an ample class on $S$. Suppose that $q$ is large so that $q^{\ell(\alpha)/(-K_S.\alpha)}>C$ and $q > C_3$ hold.  Then we have
\[
\lim_{m \to \infty}\frac{\#\mathrm{Mor}(\mathbb P^1, S, m\alpha)}{q^{-mK_S.\alpha}}  = \tau_{-K_S}(S),
\]
where $\tau_{-K_S}(S)$ is the Tamagawa number introduced in \cite{Peyre}.
In this situation, it is given by the Euler product:
\[
\tau_{-K_S}(S)= q^2 (1-q^{-1})^{-5}\prod_{c \in |\mathbb P^1|} (1-q^{-|c|})^{5}(1 +5q^{-|c|} + q^{-2|c|} ),
\]
where $|\mathbb P^1|$ is the set of closed points on $\mathbb P^1$ and $|c| = [k(c):k]$.
\end{theo}
\begin{rema}
We make $C, C_3$, and $\ell$ explicit. See Section~\ref{subsec:quintic}.
\end{rema}

This theorem can be viewed as the all height approach of Manin's conjecture originally proposed by David Bourqui (\cite[Question 1.10]{Bou} and \cite[Definition 2.2]{Bourqui13}) and Emmanuel Peyre (\cite[Question 6.27]{Peyre21}). See also \cite[Conjectures 5.16 and 5.17]{LT26} for a more modern formulation of this approach.

Then we establish the following version of Manin's conjecture by using results for all height approach:

\begin{theo}
\label{theo:intromainII}

Let $\epsilon > 0$ be a sufficiently small rational number. We define a rational polyhedral cone $\mathrm{Nef}_1(S)_\epsilon$ by
\[
\mathrm{Nef}_1(S)_\epsilon = \{ \alpha \in \mathrm{Nef}_1(S) \, | \, \ell(\alpha)/(-K_S.\alpha) \geq \epsilon \},
\]
where $\ell(\alpha)$ is a rational piecewise linear function from the previous theorem.
There exists a uniform constant $C \geq 1$ (the same constant as in Theorem~\ref{theo:intromainI}) with the following properties:
let $S$ be a smooth split quintic del Pezzo surface over $k = \mathbb F_q$ and $\epsilon > 0$ be a small rational number. Assume that $q$ is large so that $q^\epsilon > C$ holds. Then we have
\[
N(\mathbb P^1, S, -K_S, \mathrm{Nef}_1(S)_\epsilon, d) \sim \alpha(-K_S, \mathrm{Nef}_1(S)_\epsilon)\tau_{-K_S}(S)q^dd^4,
\]
as $d \to \infty$ where $\alpha(-K_S, \mathrm{Nef}_1(S)_\epsilon)$ is the alpha constant introduced in \cite{Peyre}. See also \cite[Definition 8.2]{DLTT25} for its definition.
\end{theo}

\begin{rema}
In particular, this establishes a lower bound of correct magnitude for the counting function $N(\mathbb P^1, S, -K_S, \mathrm{Nef}_1(S), d)$. We also obtained an upper bound of correct magnitude for such a counting function. See Theorem~\ref{theo:upperbound}.
\end{rema}

\subsubsection*{Low degree del Pezzo surfaces}

\cite{DLTT25} established a lower bound of correct magnitude for split quartic del Pezzo surfaces defined over finite fields. Similar results can be obtained for lower degree del Pezzo surfaces too:

\begin{theo}
\label{theo:intromainIII}
There exists a uniform constant $C_e \geq 1$ only depending on $e = 1, 2, 3$ with the following properties:
let $S$ be a smooth split del Pezzo surface of degree $e \leq 3$ defined over $k = \mathbb F_q$ such that $q > C_e$.
Then there exists a rational polyhedral cone $\mathcal C \subset \mathrm{Nef}_1(S)$ with a positive alpha constant such that we have
\[
N(\mathbb P^1, S, -K_S, \mathcal C, d) \sim \alpha(-K_S, \mathcal C)\tau_{-K_S}(S)q^dd^{9-e},
\]
as $d \to \infty$. Here $\tau_{-K_S}(S)$ is given by
\[
\tau_{-K_S}(S) = q^2 (1-q^{-1})^{-r-2}\prod_{c \in |\mathbb P^1|} (1-q^{-|c|})^{r+2}(1 +(r+2)q^{-|c|} + q^{-2|c|} ),
\]
where $r = 8-e$.
\end{theo}

\subsubsection*{Rationality}
Our results rely on the structure theorems of the moduli space $\mathrm{Mor}(\mathbb P^1, S, \alpha)$. This has another implication to the geometry of the moduli space. Here is our result:

\begin{theo}
\label{theo:intromainIV}
Let $k$ be an arbitrary field.
Let $S$ be a split smooth del Pezzo surface of degree $4\leq e \leq 7$ defined over $k$.
Let $\alpha \in \mathrm{Nef}_1(S)_{\mathbb Z}$ be a non-zero nef class.
Then the morphism scheme
\[
\mathrm{Mor}(\mathbb P^1, S, \alpha),
\]
is rational over $k$.
\end{theo}

We should mention that \cite{FGP19} established unirationality of the moduli spaces of birational stable maps of genus $0$ for del Pezzo surfaces over an algebraically closed field, but their results do not apply to all classes when $S$ is a quartic del Pezzo surface.

\subsection{History}

We will mention results on Manin's conjecture for smooth del Pezzo surfaces. Over number fields, Manin's conjecture for del Pezzo surfaces of degree $\geq 6$ has been settled by Batyrev--Tschinkel, using the fact that those surfaces are toric (\cite{BTtoric96, BTtoric98}). Manin's conjecture for del Pezzo surfaces of degree $5$ has been extensively studied,  attested by \cite{delaBre, Browning22, BD24, BD25, HBL25}, and the problem has been solved in the split case. (We should point out that \cite{Browning22, HBL25} used conic bundle structures like what we did in this paper. Other papers are based on universal torsors.)
There is only a single example of a degree $4$ del Pezzo surface where Manin's conjecture has been proved.  (\cite{DBB11}) Upper bounds and lower bounds of correct magnitude for quartic del Pezzo surfaces with conic bundle structures have been obtained in \cite{BS19}, and lower bounds of correct magnitude for del Pezzo surfaces with conic bundle structures are also available by \cite{FLS}.

Over global function fields, Glas and his collaborator obtained upper bounds for del Pezzo surfaces of lower degree in \cite{GH24, Glas25}, and \cite{DLTT25} established a weaker version of Manin's conjecture for rational curves on split quartic del Pezzo surfaces. The current paper is heavily based on results in \cite{DLTT25}.

\bigskip

\noindent
{\bf Acknowledgements:}
The author would like to thank his collaborators Ronno Das, Brian Lehmann, and Philip Tosteson for fruitful collaborations.
The author would like to thank Kazuhiro Fujiwara and Will Sawin for helpful discussions.
The author would like to thank Tim Browning for his comments on a draft of this paper.
Finally the author would like to thank the referee for helpful comments which significantly improved the exposition of the paper.

The author was partially supported by JST FOREST program Grant number JPMJFR212Z and by JSPS KAKENHI Grant-in-Aid (B) 23K25764.

%\section{Preliminaries}
%\label{sec:preliminaries}

\section{Del Pezzo surfaces}

\subsection{Notation}

Our ground field is denoted by $k$. 
A variety defined over $k$ is an integral separated scheme of finite type over $k$.

Let $X$ be a smooth projective variety defined over $k$.
We denote the space of $\mathbb R$-divisors up to numerical equivalence by $N^1(X)$ and the space of real $1$-cycles up to numerical equivalence by $N_1(X)$.
The lattice in $N_1(X)$ generated by integral classes is denoted by $N_1(X)_{\mathbb Z}$. 
We let $\mathrm{Nef}_1(X) \subset N_1(X)$ be the nef cone of $1$-cycles.
We denote the base change of $X$ to the algebraic closure by $X_{\overline{k}}$.

Let $L$ be an ample divisor on $X$. An $L$-line is a rational curve $C \subset X$ such that $L.C = 1$ An $L$-conic is a rational curve $C \subset X$ such that $L.C = 2$.

\subsection{Del Pezzo surfaces}

Here we collect several auxiliary results regarding del Pezzo surfaces.
The results in this section are based on \cite[Section 3]{DLTT25}.
We include proofs for completeness.
Let $k$ be a field and let $S$ be a smooth split del Pezzo surface of degree $e = (-K_S)^2$ defined over $k$. 
Here split means that we have $\mathrm{Pic}(S) \cong \mathrm{Pic}(S_{\overline{k}})$.
%Let $F$ be the numerical class of a smooth conic on $S$; there is an associated conic bundle structure $\pi_F : S \to \mathbb P = \mathbb P^1$.  Letting $r = 8-e$, the map $\pi_{F}$ has $r$ singular fibers. %which we denote by
%\[
%E_{1, 1} + E_{1, 2}, \cdots, E_{r, 1} + E_{r, 2},
%\]
%with images $p_1, \cdots, p_r \in \mathbb P$ respectively.  
%For any class $\alpha \in N_1(S)_{\mathbb Z}$, we define the following intersection numbers:
%\[
%h(\alpha) = -K_S.\alpha, \quad \qquad a(\alpha) = F.\alpha, \quad \qquad k_i(\alpha) = E_{i, 1}.\alpha,
%\]
%for $i = 1, \cdots, r$.   

A version of the following result is obtained for quartic del Pezzo surfaces in \cite[Section 3]{DLTT25}. Here we provide a slightly different proof:

\begin{lemm}
\label{lemm:identifyingconicfibration}
Let $S$ be a split del Pezzo surface of degree $5 \leq e \leq 7$.
Let $\alpha \in \Nef_1(S) \cap N_1(S)_{\bQ}$ be an ample $\mathbb Q$-class.

Then there exist a birational morphism $\phi : S \to S_0 \cong \mathbb P^1 \times \mathbb P^1$ and a conic fibration $\pi_F : S \to S_0 \to \mathbb P^1$ such that if we denote exceptional divisors for $\phi$ by $E_{1, 1}, \cdots, E_{r, 1}$ and a general fiber of $\pi_F$ by F, then we have
\[
-K_S.\alpha > 2F.\alpha \,  \text{ and } \, 2F.\alpha > \sum_{i = 1}^r E_{i, 1}.\alpha. 
\]
\end{lemm}

\begin{proof}
We define 
\[
h(\alpha) = -K_S.\alpha, \quad \qquad a(\alpha) = F.\alpha, \quad \qquad k_i(\alpha) = E_{i, 1}.\alpha.
\]
First let us assume that $S$ is a smooth del Pezzo surface of degree $7$.
Then $S$ admits a unique birational morphism $\phi : S \to S_0$.
Let $E_{1, 1}$ be the exceptional divisor for $\phi$. Let $F_1, F_2$ be smooth conics of two conic fibrations coming from $S_0$. Then we have $-K_S > 2F_i > E_{1, 1}$. Here for two divisors $D_1, D_2$, $D_1 > D_2$ means that $D_1-D_2$ is linearly equivalent to a non-zero effective divisor. Indeed, $-K_S - 2F_i$ and $2F_i - E_{1, 1}$ are linearly equivalent to a conic + a $(-1)$-curve.
Thus our assertion follows for both conics.

Next let us assume that $S$ is a smooth del Pezzo surface of degree $6$. We show that for any nef class $\alpha$ on $S$,  there exist a birational morphism $\phi : S \to S_0 \cong \mathbb P^1 \times \mathbb P^1$ and a conic fibration $\pi_F : S \to S_0 \to \mathbb P^1$
such that we have 
\[
h(\alpha) \geq 2a(\alpha), \, 2a(\alpha) \geq k_1(\alpha) + k_2(\alpha).
\]
Indeed there exists a blow down $\phi : S \to S'$ to a degree $7$ del Pezzo surface such that $\alpha$ can be expressed as
\[
-sK_S + \phi^*\beta,
\]
where $s$ is a non-negative rational number and $\beta$ is a nef $\mathbb Q$-class on $S'$.
Then we find a birational morphism $\phi' : S' \to S_0 \cong \mathbb P^1 \times \mathbb P^1$ and a conic fibration $\pi_F' : S' \to S_0 \to \mathbb P^1$ such that $h(\beta) \geq 2a(\beta) \geq k_1(\beta)$.
Then with respect to $\phi^*F$ we have
\[
h(\alpha) = 6s + h(\beta), \, a(\alpha) = 2s + a(\beta), k_1(\alpha) + k_2(\alpha) = 2s + k_1(\beta).
\]
Thus our assertion follows.

Next assume that $S$ is a del Pezzo surface of degree $5$. We write our ample class $\alpha$ as
\[
\alpha = -sK_S + \phi^*\beta,
\]
where $s$ is a positive rational number, $\phi : S \to S'$ be a blow down to a degree $6$ del Pezzo surface, and $\beta$ is a nef $\mathbb Q$-class on $S'$. Then we can find a birational morphism $\phi' : S' \to S_0 \cong \mathbb P^1 \times \mathbb P^1$ and a conic fibration $\pi_F' : S' \to S_0 \to \mathbb P^1$ such that $h(\beta) \geq 2a(\beta)$.
Then with respect to $\phi^*F$, we have
\[
h(\alpha) = 5s + h(\beta), \, a(\alpha) = 2s + a(\beta),  \, \sum_{i = 1}^3k_i(\alpha) = 3s + k_1(\beta) + k_2(\beta).
\]
Since $s$ is positive, we have $h(\alpha) > 2a(\alpha), 2a(\alpha) > \sum_{i = 1}^3k_i(\alpha)$.

\end{proof}

\subsection{Geometry of conic bundles}
\label{sec:geometryofconicbundles}

Let $k$ be an algebraically closed field. Here we collect some geometric results about conic bundle structures on del Pezzo surfaces. Let $S$ be a del Pezzo surface of degree $e$ defined over $k$. We assume that $S$ admits a conic fibration $\pi_F : S \to \mathbb P$ where $F$ is the class of a general fiber of $\pi_F$ and $\mathbb P\cong \mathbb P^1$.  Setting $r = 8-e$ we see that $\pi_F$ has $r$ singular fibers. We consider $\mathcal X = S\times B$ where $B = \mathbb P^1$ and $\rho_F : \mathcal X \to \mathbb P \times B$ given by $\rho_F = \pi_F \times \mathrm{id}_B$. We also denote the projection $\mathcal X \to B$ by $\pi$ and the other projection $\mathcal X \to S$ by $p$. Any rational curve $f : B \to S$ can be identified with a section $\widetilde{f} : B \to \mathcal X$ of $\pi$, and we will abuse  notation by identifying them.

Let $f : B \to \mathcal X$ be a section of $\pi$ and we denote its image by $C$. We define the height of $f$ as $h(f) = -p^*K_S.C$.
Let $p_1, \cdots, p_r \in \mathbb P$ be the images of the singular fibers of $\pi_F$.  We denote the two irreducible components of the singular fiber over $p_{i}$ respectively by
\[
E_{1, 1} + E_{1, 2}, \cdots, E_{r, 1} + E_{r, 2}.
\]
Let $B_i = \{p_i\} \times B \subset \mathbb P\times B$ for $i = 1, \cdots, r$.
We define $a(f) = p^*F.C$. Then $C' =( \rho_F)_*C$ is a section of type $(1, a(f))$ on $\mathbb P\times B$.
We denote the pullback $\rho_F^*C'$ by $Y_{C'}$. 
First we record the following lemma:

\begin{lemm}
\label{lemm: smoothness}
Let $a > 0$ be a positive integer and $C'$ be a section of type $(1, a)$ on $\mathbb P\times B$.
The surface $\rho_F|_{Y_{C'}} : Y_{C'}\to C'$ is a conic bundle with at most $ra$ singular fibers; in particular $-K_{Y_{C'}}$ is relatively ample.  When $C'$ meets with $B_i$ transversally for every $i$, $Y_{C'}$ is smooth and $\rho_F|_{Y_{C'}} : Y_{C'}\to C'$  comes with exactly $ra$ singular fibers.  In general $Y_{C'}$ admits only canonical singularities of type $A_n$.
\end{lemm}

\begin{proof}
It is clear that $Y_{C'}$ is smooth over the points $b \in C'$ which do not map to $\{p_1, \cdots, p_r\} \subset \mathbb P$.  Suppose that $b$ maps to $p_i$ so that the fiber of $\rho_F|_{Y_{C'}}$ over $b$ is reducible.  The local equation of $Y_{C'}$ near the node of this fiber has the form $xy = z^{n}$ where $n$ is the intersection multiplicity of $C'\cap B_i$ at $b$. This implies that $Y_{C'}$ has only canonical singularities and $-K_{Y_{C'}}$ is relatively ample because it intersects positively against any vertical curve. When $C'$ meets $B_i$ transversally, we have $n = 1$. Thus $Y_{C'}$ is smooth.
\end{proof}

%Let $C'$ be a section of type $(1, a)$ on $\mathbb P\times B$ with $a >0$. We denote $\rho_F^*C'$ by $Y_{C'}$. 
Let $\beta : \widetilde{Y}_{C'} \to Y_{C'}$ be the minimal resolution of singularities. Since $Y_{C'}$ has canonical singularities, by the definition, we have
\[
K_{\widetilde{Y}_{C'}} = \beta^*K_{Y_{C'}}.
\]
Let $x \in C' \cap B_i$ and suppose that the intersection multiplicity of $C'\cap B_i$ at $x$ is given by $n+ 1$. We denote the singular fiber at $x$ in $Y_{C'}$ by $E_{i, 1, x} + E_{i, 2, x}$. Then at the node, $Y_{C'}$ has a singularity of type $A_n$ so that the exceptional divisor of the minimal resolution $\beta$ at the node is given by a chain of $(-2)$-curves, $N_{1, x} + \cdots  + N_{n, x}$ such that $N_{i, x}.N_{i + 1, x} = 1$. Then the strict transform $\widetilde{E}_{i, 1, x}$ of $E_{i, 1, x}$ is attached to $N_{1, x}$ and the strict transform $\widetilde{E}_{i, 2, x}$ of $E_{i, 2, x}$ is attached to $N_{n, x}$. Thus the fiber of $\rho_F|_{Y_{C'}} : \widetilde{Y}_{C'} \to C'$ at $x$ is a chain of $n+2$ smooth rational curves.

Let $C$ be a section of $\rho_F|_{Y_{C'}} : Y_{C'} \to C'$. The height of $C$ in $Y_{C'}$ is given by
\[
h(C) = -K_{Y_{C'}/C'}.C = -K_{Y_{C'}}.C -2.
\]
Note that $C'$ is isomorphic to $\mathbb P^1$, so the degree of the tangent bundle of $C'$ is $2$.
Since $C$ is a section, $C$ satisfies that $C.p^*F = 1$.
Let $x \in C' \cap B_i$ and denote the intersection multiplicity of $C'\cap B_i$ at $x$ by $(C'\cap B_i)_x$. Then for any $i = 1, \cdots, r$ and $x \in C' \cap B_i$, $C$ satisfies
\[
 (C'\cap B_i)_xE_{i, 1, x}. C = \ell_{i, x},
\]
such that $\ell_{i, x}$ is an integer between $0$ and $(C'\cap B_i)_x$.

\begin{defi}
\label{defi:virtualclassofsections}
Let $D$ be a linear equivalence class of integral Weil divisors in $Y_{C'}$. The class $D$ is called a virtual class of sections if $D$ satisfies
\[
p^*F.D = 1, \quad  (C'\cap B_i)_yE_{i, 1, y}. D = \ell_{i, y},
\]
for all $i = 1, \cdots, r$ and $y \in C'\cap B_i$ where $\ell_{i, y}$ is an integer between $0$ and $(C'\cap B_i)_y$.
\end{defi}

Let $C$ be a section of $\rho_F|_{Y_{C'}} : Y_{C'} \to C'$ such that for any $i = 1, \cdots, r$ and $y \in C' \cap B_i$, $C$ satisfies
\[
 (C'\cap B_i)_yE_{i, 1, y}. C = \ell_{i, y}.
\]
Let $\widetilde{C}$ be the strict transform of $C$. Then $\widetilde{C}$ satisfies
\[
-K_{\widetilde{Y}_{C'}/C'}.\widetilde{C} = -K_{Y_{C'}/C'}.C, \quad \beta^*p^*F.\widetilde{C} = 1,
\]
and for all $i = 1, \cdots, r$ and $y \in C' \cap B_i$, $\widetilde{C}$ intersects with $N_{j, y}$ but not other components in the singular fiber of $\widetilde{Y}_{C'} \to C'$ at $y$ where $j = (C'\cap B_i)_y - \ell_{i, y}$. (Here $N_{0, y}$ is interpreted as $\widetilde{E}_{i, 1, y}$ and $N_{(C' \cap B_i)_y, y}$ is interpreted as $\widetilde{E}_{i, 2, y}$.) Motivated by this we introduce the following definition:

\begin{defi}
\label{defi:virtualstricttransform}
Let $D$ be a virtual class of sections for $Y_{C'}$ as above.
%such that  for any $i = 1, \cdots, r$ and $x \in C' \cap B_i$, $D$ satisfies
%\[
 %(C'\cap B_i)_xE_{i, 1, x}. D = \ell_{i, x},
%\]
%where $\ell_{i, x}$ is an integer between $0$ and $(C'\cap B_i)_x$.
Then the virtual strict transform $\widetilde{D}$ of $D$ is uniquely defined by
\[
-K_{\widetilde{Y}_{C'}/C'}.\widetilde{D} = -K_{Y_{C'}/C'}.D, \quad \beta^*p^*F.\widetilde{D} = 1,
\]
and for any $i = 1, \cdots, r$ and $y \in C' \cap B_i$, $\widetilde{D}$ intersects with $N_{j, y}$ with multiplicity $1$ but not with other components in the singular fiber of $\widetilde{Y}_{C'} \to C'$ at $y$ where $j = (C'\cap B_i)_y - \ell_{i, y}$. 

Note that it follows from this definition that $\beta_*\widetilde{D} = D$ holds.
\end{defi}

Let $D$ be a virtual class of sections on $Y_{C'}$ and $\widetilde{D}$ be the virtual strict transform of $D$. Then we have a natural pushforward map which defines a morphism between linear systems of corresponding divisorial sheaves
\[
\beta_* : |\widetilde{D}| \to |D|, D' \mapsto \beta_*D'.
\]
We claim that this is an isomorphism:
\begin{lemm}
\label{lemm:virtualstricttransform}
Let $C'$ be a section of type $(1, a)$ on $\mathbb P \times B$ with $a > 0$. Let $D$ be a virtual class of sections on $Y_{C'}$ and $\widetilde{D}$ be the virtual strict transform of $D$. Then the pushforward map
\[
\beta_* : |\widetilde{D}| \to |D|, D' \mapsto \beta_*D'.
\]
is an isomorphism.
\end{lemm}
\begin{proof}
We show that for any $D' \in |D|$, there exists a unique $\widetilde{D}' \in |\widetilde{D}|$ such that $\beta_*\widetilde{D}' = D'$. First let us prove uniqueness. Let $\widetilde{D}_1', \widetilde{D}_2' \in |\widetilde{D}|$ be two effective divisors on $\widetilde{Y}_{C'}$ such that $\beta_*\widetilde{D}_1' = \beta_*\widetilde{D}_2'$. Then $D_1' - D_2'$ is supported on exceptional divisors of $\beta$. Since $\widetilde{D}_1'$ and $\widetilde{D}_2'$ have the same intersection properties, it follows from \cite[Lemma 1-2-10]{Matsuki02} that $\widetilde{D}_1' = \widetilde{D}_2'$.

Next we prove the existence. If $D'$ is an irreducible section, then we take $\widetilde{D}'$ as the strict transform of $D'$. Suppose that $D'$ is reducible so that it is the sum of a section $C$ and an effective  vertical divisor. For a full fiber, one can take the corresponding full fiber on $\widetilde{Y}_{C'}$ as the strict transform. Thus let us demonstrate this claim when $D' = C + tE_{i, 1, x}$ for some $i$ and $x \in C'\cap B_i$. Let $n + 1 = (C'\cap B_i)_x$, $\ell = (n+1)D'.E_{i, 1, x}$, and $m = (n+1)C.E_{1, i, x}$.
Then since we have $E_{i, 1, x}^2 = -\frac{1}{n+1}$, we conclude $t = m-\ell$. Let $\widetilde{C}$ be the strict transform of $C$. Then 
\[
\widetilde{D}' = \widetilde{C} + \sum_{j = n-m}^{n-\ell}N_{j, x},
\]
witnesses our claim. Thus our assertion.
\end{proof}

\section{The space of rational curves on del Pezzo surfaces}
\label{sec:geometryofmoduli}

In this section we study some basic properties of the space of rational curves on del Pezzo surfaces.
We work over a field $k$ of arbitrary characteristic.

\subsection{The transversality of the space of sections/the space of rational curves}

Let $S$ be a smooth del Pezzo surface of degree $e$ defined over $k$.
%We will denote the Picard rank of $S$ by $\rho(S)$ and set $r = \rho(\overline{S}) -2$.
We set $\pi : \mathcal X = S \times B \to B$ where $B = \mathbb P^1$.
Let $p : \mathcal X \to S$ be the other projection.
We also denote any fiber of $\pi$ by $S_b$.

We denote {\it the space of sections of $\pi$} by $\mathrm{Sec}(\mathcal X)$. This can be identified with the space of rational curves $\mathrm{Mor}(B, S)$ via the composition $$(f : B \to\mathcal X )\mapsto (p\circ f : B \to S).$$
We consider the nef cone of curves $\mathrm{Nef}_1(\mathcal X)$ and define $$\mathrm{Nef}_1(\mathcal X)_{S = 1}= \{\alpha \in \mathrm{Nef}_1(\mathcal X) \, |\, \alpha.S_b = 1\}.$$
For each $\alpha \in \mathrm{Nef}_1(\mathcal X)_{S = 1}$ we define the space of sections of class $\alpha$ by $\Sec(\mathcal X/B, \alpha)$. 
We also denote the height of $\alpha$ by 
\[
h(\alpha):= -p^*K_S.\alpha.
\]
Then we have the following theorem:

\begin{theo}[{\cite{Testa05}, \cite{Testa09}, and \cite{BLRT21}}]
\label{theo: Testa}
Let $$\alpha \in \mathrm{Nef}_1(\mathcal X)_{S = 1} \cap N_1(\mathcal X)_{\mathbb Z},$$ such that $-p^*K_S.\alpha \geq 2$. When the characteristic $p$ of $k$ is positive, we furthermore assume that
\[
p \geq
\begin{cases}
2 & \text{if $e \geq 4$;} \\
3 & \text{if $e = 3, 2$, or;}\\
11 & \text{if $e = 1$.}
\end{cases}
\]

Assume that $p_*\alpha$ is not a multiple of the class of a $-K_S$-conic on $S$. Then the space $\Sec(\mathcal X/B, \alpha)$ consists of a unique dominant component which has the expected dimension 
\[
-K_{\mathcal X/B}.C + 2 = -p^*K_S.C + 2 = h(\alpha) + 2,
\]
i.e., the dimension of this unique component is equal to the expected dimension.

Next assume that $p_*\alpha$ is a multiple of the class of a $-K_S$-conic $C$ on $S$. 
If $C^2 = 0$ then the space $\Sec(\mathcal X/B, \alpha)$ consists of a unique dominant component of expected dimension. 

If $C^2 = 2$ (this can only happen when $e = 1$ or $2$), then the space $\Sec(\mathcal X/B, \alpha)$ consists of exactly two dominant components of expected dimension. 

If $C^2 = 4$ (this can only happen when $e = 1$) and $S$ is general, then the space $\Sec(\mathcal X/B, \alpha)$ consists of exactly two dominant components of expected dimension. 
Even when $S$ is not general, the space $\Sec(\mathcal X/B, \alpha)$ consists of at least two dominant components and every component has expected dimension.
\end{theo}

A version of the above statement for higher genus curves is obtained in \cite{LT21} in characteristic $0$. 
\cite{LT19} proved the irreducibility of the space of sections for certain non-isotrivial del Pezzo fibrations over $\mathbb P^1$.

\subsection{Fibrations structures induced by conic bundle structures}

Next we discuss some structural results on $\Sec(\mathcal X/B, \alpha)$.
Let $F$ be a class of smooth conics on $S$.
We denote the linear system of $F$ by $\mathbb P$ which is isomorphic to $\mathbb P^1$.
We consider the conic fibration $\pi_F : S \to \mathbb P$ induced by $F$ and this induces $\rho_F = \pi_F \times \mathrm{id}_B : \mathcal X \to \mathbb P \times B$.

We fix a class $\alpha \in \mathrm{Nef}_1(\mathcal X)_{S = 1} \cap N_1(\mathcal X)_{\mathbb Z}$. We denote $a(\alpha) = p^*F.\alpha$. For each $C \in \Sec(\mathcal X/B, \alpha)$ the pushforward $(\rho_F)_*C$ is a section of type $(1, a(\alpha))$ on $\mathbb P\times B$.
We denote the space of sections of type $(1, a(\alpha))$ on $\mathbb P\times B$ by $\mathrm{Sec}(\mathbb P\times B/B, a(\alpha))$.
\begin{lemm}
\label{lemm: affineness}
Assume that $a(\alpha) > 0$.
The space $\mathrm{Sec}(\mathbb P\times B/B, a(\alpha))$ is a Zariski open set of the projective space $|(\rho_F)_*\alpha|$ of dimension $2a(\alpha) + 1$. Its complement is irreducible of codimension $1$ and generically parametrizes unions of a section of type $(1, a(\alpha)-1)$ and a fiber $P$.
In particular $\mathrm{Sec}(\mathbb P\times B/B, a(\alpha))$ is an affine variety.
\end{lemm}

\begin{proof}
We may assume that $k$ is algebraically closed.
All sections of type $(1, a(\alpha))$ on $\mathbb P \times B$ are linearly equivalent to each other so they form a Zariski open subset $\mathrm{Sec}(\mathbb P\times B/B, a(\alpha))$ of the linear system $\mathbb P_{a(\alpha)} : = |C'| = \mathbb P^{2a(\alpha) + 1}$ where $C'$ is a section of type $(1, a(\alpha))$.
We denote the complement of $\mathrm{Sec}(\mathbb P\times B/B, a(\alpha))$ in $\mathbb P_{a(\alpha)}$ by $\mathcal D$.
Then the morphism
\[
\mathbb P_{a(\alpha)-1} \times B \to \mathcal D, (C'', b) \mapsto C''+ \mathbb P_b,
\]
where $\mathbb P_b$ is the fiber of the projection $\mathbb P\times B \to B$ at $b$, is the normalization of $\mathcal D$. Thus $\mathcal D$ is irreducible of codimension $1$, proving the claim.
\end{proof}

We have a natural morphism
\[
\Phi_\alpha : \Sec(\mathcal X/B, \alpha) \to \mathrm{Sec}(\mathbb P\times B/B, a(\alpha)), \quad C \mapsto (\rho_F)_*C.
\]
Our goal is to understand the structures of this morphism. 
Let $$C' \in \mathrm{Sec}(\mathbb P\times B/B, a),$$ and we denote the surface $\rho_F^*C'$ by $Y_{C'}$. The linear system $\mathbb P$ admits $r(S)$ points representing singular conics and we denote them by $p_1, \cdots, p_{r(S)}$. We set $B_i = \{p_i\} \times B \subset \mathbb P\times B$ for every $i = 1, \cdots, r(S)$.

Next we discuss when $\Phi_\alpha$ is dominant for $\alpha \in \mathrm{Nef}_1(\mathcal X)_{S = 1} \cap N_1(\mathcal X)_{\mathbb Z}$:
\begin{prop}
\label{prop:dominant}
When the characteristic $p$ of $k$ is positive, we assume that either $e \geq 2$ or $e= 1$ and $p \geq 11$. Furthermore we may assume that our del Pezzo surface $S$ is not isomorphic to one of surfaces specified in \cite[Theorem 1.2]{BLRT21}.

Let $\alpha \in \Nef_1(\mathcal X)_{S=1} \cap N_1(\mathcal X)_{\bZ}$ be a nef class such that $h(\alpha) = -p^*K_S.\alpha \geq 2$. Let $a(\alpha) = p^*F.\alpha$. Then
the map $\Phi_\alpha$ is dominant when $h(\alpha) \geq 2a(\alpha)$.
\end{prop}

\begin{proof}

The statement is purely geometric, so we may assume that $k$ is algebraically closed after taking a base change to an algebraic closure.
It follows from \cite[Theorem 1.2]{BLRT21} that $\Sec(\mathcal X/B, \alpha)$ is non-empty and each component parametrizes a free rational curve. In particular a general $C \in \Sec(\mathcal X/B, \alpha)$ is a free curve on $\mathcal X$ and the dimension of each component is given by $h(\alpha) + 2$.

For a general $C \in \Sec(\mathcal X/B, \alpha)$, we denote its image $(\rho_F)_*C$ by $C'$ and the surface $\rho_F^*C'$ by $Y_{C'}$. Let $\widetilde{Y}_{C'}$ be the minimal resolution of $Y_{C'}$ 
and let $\widetilde{C}$ be the strict transform of $C$ in $\widetilde{Y}_{C'}$.
Then the degree of the normal bundle of $\widetilde{C}$ as a rational curve in $\widetilde{Y}_{C'}$ is $h(\alpha) -2a(\alpha) \geq 0$. Indeed, we have
\[
\deg(N_{\widetilde{C}/\widetilde{Y}_{C'}}) = -K_{\widetilde{Y}_{C'}/C'}.\widetilde{C} = -K_{Y_{C'}/C'}.C = -p^*K_S.C - Y_{C'}.C = h(\alpha) -2a(\alpha).
\]
Thus it is free in $\widetilde{Y}_{C'}$, and the expected dimension and the actual dimension of this moduli space agree and they are equal to $h(\alpha) -2a(\alpha) + 1$.  We conclude that every component of $\Sec(\mathcal X/B, \alpha)$ maps to $\Sec(\mathbb P\times B/B, a)$ dominantly: the dimension of $\Sec(\mathcal X/B, \alpha)$ is $h(\alpha) + 2$, the dimension of $\Sec(\mathbb P\times B/B, a)$ is $2a(\alpha)+1$, and we just proved that the dimension of a general fiber of $\Phi_{\alpha}$ is the difference.
\end{proof}

\subsection{Configuration covers}
\label{subsec: configuration}
The discussion in the remaining sections is based on \cite[Section 4]{DLTT25}.
After taking a partial compactification of $\Sec(\mathcal X/B, \alpha)$ to make $\Phi_{\alpha}$ proper, the Stein factorization of $\Phi_{\alpha}$ defines a finite cover $\phi_\alpha : Z_\alpha \to \mathrm{Sec}(\mathbb P\times B/B, a(\alpha))$.
We are going to compute the degree of $\phi_\alpha$.
We denote singular fibers of $\pi_F : S \to \mathbb P$ by
\[
E_{1, 1}+ E_{1, 2}, \cdots, E_{r, 1}+ E_{r, 2}.
\]
We denote $k_i(\alpha) = \alpha.p^*E_{i, 1}$ for $i = 1, \cdots, r(S)$. 
\

To discuss the degree of $\phi_\alpha$, let us introduce a certain configuration cover.
Let $\mathbb P_a$ be the linear system of divisors of type $(1, a)$ in $\mathbb P\times B$ with $a \geq 0$.  Recall that we let $p_1, \cdots, p_{r} \in \mathbb P$ be points corresponding to singular conics of $\pi_F$
and let $B_i = \{p_i\} \times B \subset \mathbb P\times B$.
Let $V_{i, a}$ be the locus of $\mathbb P_a$ parametrizing the union of a section $B_i$ and $a$ vertical fibers. Note each $V_{i,a}$ is a linear subspace of dimension $a$ and that these subspaces are mutually disjoint. 

We fix $i = 1, \cdots, r$ and $0\leq k_i \leq a$ and we define
\[
Z_{a, i, k_i} = \{([C], [T]) \in \mathbb P_a \times \mathrm{Hilb}^{[k_i]}(B_i)\, | \, T \subset C\} \subset \mathbb P_a \times \mathrm{Hilb}^{[k_i]}(B_i).
\]

Let $\bold k = (k_1, \cdots, k_r)$.
When $h(\alpha) \geq 2a(\alpha)$, it follows from Proposition~\ref{prop:dominant} that our space $Z_\alpha$ can be considered as a Zariski open subset of the scheme
\[
Z_{a, \bold k} =  Z_{a, 1, k_1} \times_{\mathbb P_a} \cdots \times_{\mathbb P_a} Z_{a, r, k_r}\subset \mathbb P_a \times \prod_{i = 1}^{r(S)}\mathrm{Hilb}^{[k_i]}(B_i),
\]
with $a = a(\alpha)$ and the projection $Z_{a, \bold k} \to \mathbb P_a$ by $\pi_{a, \bold k}$. We call these varieties {\it configuration covers}. This is originally introduced in \cite[Definition 4.1]{DLTT25}. We denote the preimage of $\mathbb{P}_{a}^{\circ}$ by $Z_{a, \bold k}^{\circ}$.

The relation of $Z_{a, \bold k}$ to the space of sections is given by the following morphism:
We define the morphism by
\begin{align}
\label{structuremorphism}
\psi_\alpha: \Sec(\mathcal X/B, \alpha) \to Z_{a(\alpha), \bold k(\alpha)}: [C] \mapsto ([(\rho_F)_*C], ((\rho_F)_*(C\cap p^*E_{i, 1}))_{i = 1}^{r(S)}),
\end{align}
where $\bold k(\alpha) = (k_1(\alpha), \cdots, k_{r(S)}(\alpha))$.
We know that this morphism is dominant when $h(\alpha) \geq 2a(\alpha)$ by Proposition~\ref{prop:dominant}. Using this we have the following lemma:
\begin{lemm}
\label{lemm:structuremorphism}
When the characteristic $p$ of $k$ is positive, we assume that either $e \geq 2$ or $e = 1$ and $p \geq 11$. Furthermore we may assume that our del Pezzo surface $S$ is not isomorphic to one of surfaces specified in \cite[Theorem 1.2]{BLRT21}. 

Let $\alpha \in \Nef_1(\mathcal X)_{S=1} \cap N_1(\mathcal X)_{\bZ}$ be a nef class such that $h(\alpha)  \geq  \max\{a(\alpha), 2\}$. 
Then $\phi_\alpha$ is separable and the degree of $\phi_\alpha$ is equal to
\[
\prod_{i = 1}^{r(S)} \binom{a(\alpha)}{k_i(\alpha)}.
\]
\end{lemm}

\begin{proof}
We may assume that $k$ is algebraically closed.
Let $$C \in \mathrm{Sec}(\mathcal X/B, \alpha),$$ be a general section. It follows from Proposition~\ref{prop:dominant} that 
\[
N_{C/\mathcal X} \to \rho_F^*N_{C'/(\mathbb P\times B)} \to 0,
\]
is exact because for a general $C$, $C'$ meets each $B_i$ transversally so $C$ is contained in the smooth locus of $\mathcal X \to \mathbb P\times B$.
Thus the kernel $K$ of the first map is a line bundle of degree $h(\alpha) - 2a(\alpha)\geq 0$.
In particular this implies that $H^0(C, N_{C/\mathcal X}) \to H^0(C, \rho_F^*N_{C'/(\mathbb P\times B)})$ is surjective. We conclude that
\[
\mathrm{Sec}(\mathcal X/B, \alpha) \to \mathrm{Sec}(\mathbb P\times B, a(\alpha)),
\]
is generically smooth. Hence it follows that $\phi_\alpha$ is separable.
%We have a natural morphism
%\[
%\psi_\alpha: \Sec(\mathcal X/B, \alpha) \to Z_{a(\alpha), \bold k(\alpha)}: [C] \mapsto ([(\rho_F)_*C], ((\rho_F)_*(C\cap p^*E_{i, 1}))_{i = 1}^{r(S)}).
%\]
A general fiber of this morphism is a Zariski open subset of a projective space of dimension $h(\alpha) + 1-2a(\alpha)$. Indeed, this is the linear system of $C$ in $Y_{C'}$ where $C' = (\rho_F)_*C$.
Thus our assertion follows from the irreducibility of $Z_{a(\alpha), \bold k(\alpha)}^\circ$.
\end{proof}
\begin{rema}
When $h(\alpha) \geq 2a(\alpha)$ fails, $\psi_\alpha: \Sec(\mathcal X/B, \alpha) \to Z_{a(\alpha), \bold k(\alpha)}$ is injective onto its image, but not dominant in general.
\end{rema}

\subsection{Basic properties of configuration covers}

Here we establish basic properties of configuration covers.
First we prove that $Z_{a, \bold k}$ is lci of expected dimension.
In this section, we always assume $k_i \leq a$ for any $i$.

\begin{prop}[{\cite[Proposition 4.2]{DLTT25}}]
\label{prop:locallycompleteintersectionZ}

The space $Z_{a, \bold k}  \subset \mathbb P_a \times \prod_{i = 1}^{r}\mathrm{Hilb}^{[k_i]}(B_i)$ is the zero locus of a section of a globally generated vector bundle of rank $\sum_i k_i$. 
\end{prop}
\begin{proof}
\cite[Proposition 4.5]{DLTT25} is stated when $r = 4$, but the same proof proves this.
\end{proof}

\begin{coro}
The space $Z_{a, \bold k}$ is a lci reduced scheme of dimension $2a+1$.
\end{coro}
\begin{proof}
We may assume that $k$ is algebraically closed.
Proposition \ref{prop:locallycompleteintersectionZ} shows that every irreducible component of $Z_{a, \bold k}$ has dimension $\geq 2a+1$.  Note that $\pi_{a, \bold k}^{-1}(Z_{a, \bold k} \setminus \cup V_{i, a})$ has dimension exactly $2a+1$.  Its complement is the disjoint union of the various $\pi_{a, \bold k}^{-1}(V_{i,a})$; the preimage of $V_{i,a}$ has dimension
\begin{equation*}
\dim(\Hilb^{[k_{i}]}(B_{i})) + \dim(V_{i,a}) = k_{i} + a < 2a+1.
\end{equation*}
We conclude that $Z_{a, \bold k}$ is of dimension $2a+1$.  Since $Z_{a, \bold k}$ is a zero section of a vector bundle with the expected codimension, it is locally defined by a regular sequence.  In particular it is reduced if and only if it is generically reduced.  But it is easy to see that $Z_{a, \bold k}$ is generically smooth (see tangent space computation later), thus reduced everywhere.

\end{proof}

Next we analyze singularities of the configuration covers $Z_{a, \bold k}$.
In particular we show that $Z_{a, \bold k}$ is normal.

Let us assume that our ground field $k$ is algebraically closed.
Let $\widetilde{Z}$ be the Hilbert-flag scheme parametrizing $([C], [T_1], \cdots, [T_r])$ where $C\in \mathbb P_a$ and $T_i \subset C$ is a length $k_i$ subscheme. Then $Z_{a, \bold k}$ is a closed subscheme of $\widetilde{Z}$.  %We will use this inclusion to analyze the singularities of $Z_{a, \bold k}$.
%We first need to understand $\widetilde{Z}$.  
%This is an irreducible scheme of dimension $2a + 1 + \sum_{i = 1}^{r} k_i$. 
%First we will prove that this scheme is smooth at $([C], ([T_i]))$ with $C$ being reduced. 
%The tangent space for the Hilbert scheme parametrizing the curves $C$ is given by
%\[
%H^0(C, N_{C/\mathbb P\times B}) = H^0(C, \mathcal O(1, a)|_C),
%\]
%which is $2a + 1$ dimensional. 
%We have a natural map
%\[
%H^0(C,  N_{C/\mathbb P\times B}) \to  H^0(T_i, N_{C/\mathbb P\times B}|_{T_i}),
%\]
%where the latter space is $k_i$-dimensional.  This map is surjective when $C$ is a section because $k_i\leq a$. Hence we conclude that
%\[
%H^0(\mathbb P\times B, \mathcal O(1, a))\to H^0(T_i, N_{C/\mathbb P\times B}|_{T_i}),
%\]
%is surjective.  
%The tangent space for the Hilbert scheme parametrizing $T_i$ in $\mathbb P\times B$ is given by
%\[
%H^0(T_i, N_{T_i/\mathbb P \times B}),
%\]
%which is $2k_i$-dimensional. We claim that we have a surjective map
%\[
%H^0(T_i, N_{T_i/\mathbb P \times B}) \to H^0(T_i, N_{C/\mathbb P\times B}|_{T_i}).
%\]
%Indeed, let $I_{T_i}, I_C$ be the ideal sheaves for $T_i$ and $C$ respectively.
%It suffices to show that a natural map
%\[
%(I_C/I_C^2)|_{T_i} \to I_{T_i}/I_{T_i}^2,
%\]
%is injective. Since $C$ is a reduced effective divisor on $\mathbb P\times B$, we have the inclusion:
%\[
%I_C/I_C^2  \hookrightarrow \Omega^1_{\mathbb P\times B}|_C,
%\]
%and using a fact that $C$ is reduced, it is also true that 
%\[
%(I_C/I_C^2)|_{T_i}  \hookrightarrow \Omega^1_{\mathbb P\times B}|_{T_i},
%\]
%is injective. Our assertion follows from this.
First note that the tangent space of $\widetilde{Z}$ at $(C, ([T_i]))$ is given by
\[
H^0(C, N_{C/\mathbb P\times B}) \times_{\prod_{i = 1}^r H^0(T_i, N_{C/\mathbb P\times B}|_{T_i})} \prod_{i = 1}^r H^0(T_i, N_{T_i/\mathbb P \times B}).
\]
Also we have
\[
Z_{a, \bold k} = \widetilde{Z} \times_{\prod_{i = 1}^r \Hilb^{[k_i]}(\mathbb P\times B)} \prod_{i = 1}^r\Hilb^{[k_i]} (B_i).
\]
Now the tangent space of $\Hilb^{[k_i]}(B_i)$ is given by $H^0(T_i, N_{T_i/B_i})$ which can be considered as a subspace of $H^0(T_i, N_{T_i/\mathbb P \times B})$. Thus the tangent space of the above fiber product is given by 
\begin{align*}
& \left( H^0(C, N_{C/\mathbb P\times B}) \times_{\prod_{i = 1}^r H^0(T_i, N_{C/\mathbb P\times B}|_{T_i})} \prod_{i = 1}^r H^0(T_i, N_{T_i/\mathbb P \times B}) \right) \\ & \qquad \qquad \qquad \times_{\prod_{i = 1}^r H^0(T_i, N_{T_i/\mathbb P \times B})} \prod_{i = 1}^r H^0(T_i, N_{T_i/B_i}),
\end{align*}
which is isomorphic to
\begin{align}
\label{tangentspaceforconfiguration}
H^0(C, N_{C/\mathbb P\times B}) \times_{\prod_{i = 1}^r H^0(T_i, N_{C/\mathbb P\times B}|_{T_i})} \prod_{i = 1}^r H^0(T_i, N_{T_i/B_i}).
\end{align}
Here the projection map for the rightmost term is defined by the following compositions:
\[
\phi_{C, T_i}: H^0(T_i, N_{T_i/B_i}) \subset H^0(T_i, N_{T_i/\mathbb P \times B}) \to H^0(T_i, N_{C/\mathbb P\times B}|_{T_i}).
\]
Let us try to understand these maps in more detail.  Let $I_C$, $I_{T_i}$, $I_{B_i}$ be the ideal sheaves for $C$, $T_i$, and $B_i$ respectively.  We assume that $C$ is not contained $V_{i, a}$ for any $i$. (In this case, the above map is just the zero map.)
Then the dual of $\phi_{C, T_i}$ is given by
\[
I_C/I_C^2 \otimes O_{T_i} \to I_{T_i}/(I_{T_i}^2 + I_{B_i}).
\]
Both spaces are $k_i$-dimensional.
So the cokernel $K$ is given by
\[
I_{T_i}/(I_{T_i}^2 + I_C + I_{B_i}) = I_{T_i}/(I_{T_i}^2 + I_{C\cap B_i}).
\]
Write $$C\cap B_i = \sum_{j= 1}^s n_j p_j,$$ 
with $n_j > 0$ and $\sum_{j = 1}^s n_j = a$, and $$T_i = \sum_{j = 1}^t m_j p_j,$$ with $t \leq s$, $0< m_j \leq n_j$ and $\sum_{j = 1}^t m_j = k_i$. Then the dimension of the cokernel is given by $$\mathrm{dr}(C, T_i) := \sum_{j = 1}^t \min\{ 2m_j, n_j\} - m_j = \sum_{j = 1}^t \min\{2m_j, n_j\} - k_i.$$ We call this number the {\it dropping rank}. Note that this always satisfy $\mathrm{dr}(C, T_i)  \leq \min\{k_i, a-k_i\}$.

Returning to $Z_{a, \bold k}$, note that a point $(C,T_{1},\ldots,T_{r})$ is smooth if and only if the image of $H^{0}(C,N_{C/{\mathbb P \times B}})$ and the image of $\prod_{i = 1}^r H^0(T_i, N_{T_i/B_i})$ meet transversally in $\prod_{i = 1}^r H^0(T_i, N_{C/\mathbb P\times B}|_{T_i})$.  In other words, the vector space $$H^{0}(C,N_{C/\mathbb P \times B}),$$ must surject onto the cokernel $K$ of the map $\prod_{i} \phi_{C,T_{i}}$.  In particular this implies that $\sum_{i} \mathrm{dr}(C, T_i) \leq 2a+1$.

Let $Q$ be the subscheme of $\cup_{i} C \cap B_{i}$ such that for each point $p_j \in \cup_{i} C \cap B_{i}$ the length of $Q$ at $p_j$ is given by $\min\{ 2m_j, n_j\} - m_j$.  A necessary and sufficient condition for smoothness is that we have a surjection $H^{0}(C,N_{C/\mathbb P \times B}) \to H^{0}(Q,N_{C/\mathbb P \times B}|_{Q})$.

Even when the rank of $\phi_{C, T_i}$ is not full rank, if these dropping ranks impose independent conditions on $H^0(C, N_{C/\mathbb P\times B})$, then the tangent space has expected dimension. 
Moreover note that when $e\geq 4$, we may assume that the sum of dropping ranks never exceed $2a$.
Thus we can conclude the following propositions: 

\begin{prop} 
\label{prop:Z_singularities}
Assume that $k$ is algebraically closed and $e \geq 4$.
Let $(C, (T_i)) \in Z_{a, \bold k}$ be a point. 
Assume that $C$ is not contained in $V_{i, a}$ for any $i$. Then
a point $(C, (T_i)) \in Z_{a, \bold k}$ is a singular point if and only if 
\begin{itemize}
\item $C$ admits a non-reduced vertical fiber $n\mathbb P_1$ of length $n > 1$, and;
\item there are at least three $i$'s such that $T_i$ intersects with $\mathbb P_1$ and the length of $T_i$ at that point is strictly less than the multiplicity of $C\cap B_i$ at that point. 
\end{itemize}

Next assume that $C$ is contained $V_{i, a}$ for some $i$, i.e., $C$ is a sum of $B_i$ with $a$ vertical fibers.
Then
a point $(C, (T_i)) \in Z_{a, \bold k}$ is a singular point if and only if 
\begin{itemize}
\item $C$ admits a non-reduced vertical fiber $n\mathbb P_1$ of length $n > 1$, and;
\item  there are at least three $j$'s such that $T_j$ intersects with $P_1$ and when $j \neq i$, the length of $T_j$ at that point is strictly less than $n$. 
\end{itemize}
\end{prop}

\begin{prop} 
\label{prop:Z_singularitiesII}

Assume that $k$ is algebraically closed and $e \leq 3$.
Let $(C, (T_i)) \in Z_{a, \bold k}$ be a point. 
Assume that $C$ is not a sum of $B_i$ with $a$ vertical fibers for any $i$. 
Then a point $(C, (T_i)) \in Z_{a, \bold k}$ is a singular point if and only if 
\begin{itemize}
\item $C$ admits a non-reduced vertical fiber $n\mathbb P_1$ of length $n > 1$ and there are at least three $i$'s such that $T_i$ intersects with $\mathbb P_1$ and the length of $T_i$ at that point is strictly less than the multiplicity of $C\cap B_i$ at that point, or;
\item $\sum_{i = 1}^r \mathrm{dr}(C, T_i) > 2a + 1$.
\end{itemize}

Next assume that $C$ is a sum of $B_i$ with $a$ vertical fibers.
Then a point $(C, (T_i)) \in Z_{a, \bold k}$ is a singular point if and only if $C$ admits a non-reduced vertical fiber $n\mathbb P_1$ of length $n > 1$ and there are at least three $j$'s such that $T_j$ intersects with $\mathbb P_1$ and when $j \neq i$, the length of $T_j$ at that point is strictly less than $n$. 

\end{prop}

This singularity analysis will lead to the following corollary:

\begin{coro}
\label{coro:normal}

The space $Z_{a, \bold k}$ is normal.
\end{coro}
\begin{proof}
We may assume that $k$ is algebraically closed.
When $e \geq 4$ a point $(C, (T_i)) \in Z_{a, \bold k}$ is singular only if $C$ admits a non-reduced vertical curve. This is a codimension $3$ condition. Finally $Z_{a, \bold k}$ is lci so our assertion follows.

When $e \leq 3$, we consider a point $(C, (T_i)) \in Z_{a, \bold k}$. We only need to show that $Z_{a, \bold k}$ is smooth in codimension $1$ so we may assume that $C$ is a section or a section + a vertical fiber. Then $(C, (T_i)) \in Z_{a, \bold k}$ is singular only if $\sum_{i = 1}^r \mathrm{dr}(C, T_i) > 2a + 1$. This can only happen when $C$ is tangent to five $B_i$'s. This is codimension $5$ condition. Thus our assertion follows.
\end{proof}

\subsection{Further structures of $Z^{\circ}_{a, \bold k}$}
\label{subsec:furtherZ}

Recall that our variety $Z_{a, \bold k}$ is embedded as
\[
Z_{a, \bold k} \subset \mathbb P_a \times \prod_{i = 1}^r \mathrm{Hilb}^{[k_i]}(B_i).
\]
Then we consider the following fibration:
\[
\phi_{a, \bold k} : Z_{a, \bold k} \to  \prod_{i = 1}^r \mathrm{Hilb}^{[k_i]}(B_i).
\]
\begin{lemm}
\label{lemm:phidominant}

Assume that $2a \geq \sum_i k_i$. Then $\phi_{a, \bold k}$ is dominant.
Moreover every fiber is a projective space of dimension $\geq 2a + 1 - \sum_{i = 1}^r k_i$.
\end{lemm}
\begin{proof}
We may assume that $k$ is algebraically closed.
For any $(T_i)\in  \prod_{i = 1}^r \mathrm{Hilb}^{[k_i]}(B_i)$, the set
\[
\{ [C] \in \mathbb P_a \, | \, \text{for any $i = 1, \cdots, r$, $C\cap B_i$ contains $T_i$ as a closed subscheme.} \}
\]
has dimension $\geq 2a + 1 - \sum_{i = 1}^r k_i >0$. Thus $\phi_{a, \bold k}$ is surjective.
Moreover every fiber is a linear subspace of $\mathbb P_a$, proving the claim.
\end{proof}

Combining above, we obtain:

\begin{coro}
The configuration cover $Z_{a, \bold k}$ is irreducible assuming $2a \geq \sum_i k_i$.
\end{coro}
\begin{proof}
Lemma~\ref{lemm:phidominant} shows that $Z_{a, \bold k}$ is connected. Thus our assertion follows form Corollary~\ref{coro:normal}.
\end{proof}

From now on, we assume that $2a \geq \sum_i k_i$.
Note that for any $i = 1, \cdots, r$, $B_i$ is naturally identified with $B$. Using this we define the following Zariski open subset 
\[
U_{\bold k} := \{ (T_i) \in \prod_{i = 1}^r \mathrm{Hilb}^{[k_i]}(B_i) \, | \, (*) \},
\]
where $(*)$ means that for any $i \neq j$, $T_i$ and $T_j$ do not share any point in their supports.
Then we have
\begin{lemm}[{\cite[Lemma 4.4]{DLTT25}}]
\label{lemm:configurationcover_fibration}
Assume that $2a \geq \sum_i k_i$.
We have
\[
\phi_{a, \bold k}(Z^{\circ}_{a, \bold k}) \subset U_{\bold k}.
\]
Moreover the fibration
\[
\phi_{a, \bold k} : \phi_{a, \bold k}^{-1}(\phi_{a, \bold k}(Z^{\circ}_{a, \bold k}) ) \to \phi_{a, \bold k}(Z^{\circ}_{a, \bold k}),
\]
is a $\mathbb P^{2a + 1 - \sum_i k_i}$-bundle over the Zariski open subset $\phi_{a, \bold k}(Z^{\circ}_{a, \bold k})$.
\end{lemm}
\begin{proof}
This is stated in \cite[Lemma 4.4]{DLTT25} for $r = 4$, but the proof actually shows this.

\end{proof}

We also need the following lemma for a later application:

\begin{lemm}[{\cite[Lemma 4.11]{DLTT25}}]
\label{lemm:dimensionestimate}
Assume $2a + 1 -\sum_{i = 1}^r k_i>0$.
Then $U_{\bold k} \setminus \phi_{a, \bold k}(Z^{\circ}_{a, \bold k})$ is empty when $\sum_{i = 1}^rk_i \leq a$ and otherwise it has dimension at most $\sum_{i = 1}^r k_i - (2a+1 - \sum_{i = 1}^r k_i) $\end{lemm}

\subsection{Fibration structures for the space of rational curves}
\label{subsec:fibrationstr}

We work over a field $k$. 
%We recall the settings. Let $S$ be a split smooth del Pezzo surface of degree $e$ and $B = \mathbb P^1$. We let $\mathcal X = S \times B$ with the projection $\pi : \mathcal X \to B$ and $p : \mathcal X \to S$. Let $\alpha \in \Nef (\mathcal X)_{S = 1}$ be a nef class whose intersection against $S_b = \pi^{-1}(b)$ is equal to $1$. 
%We are interested in the space of sections of $\pi$:
%\[
%\Sec(\mathcal X/B, \alpha).
%\]
%In particular we are interested in cohomology of this moduli space. We denote the height by $h(\alpha) = -K_{\mathcal X/B}.\alpha$. 
We fix a blow down $\phi : S \to S_0 \cong \mathbb P^1 \times \mathbb P^1$.
We have the morphism
\[
\psi_\alpha : \Sec(\mathcal X/B, \alpha) \to Z_{a(\alpha), \bold k(\alpha)}^{\circ},
\]
defined in (\ref{structuremorphism}). According to Proposition~\ref{prop:dominant}, this morphism is dominant when $h(\alpha) \geq 2a(\alpha)$. 
%Our first goal of this section is to understand the image of this morphism in general. First we discuss when $\psi_\alpha$ is surjective.

We record the following results on the dimension of the complement of the image of $\Sec(\mathcal X/B, \alpha)$ in $U_{\bold k(\alpha)}$:

\begin{prop}[{\cite[Proposition 4.13]{DLTT25}}]
\label{prop:codimensionI}
Assume that $h(\alpha) \geq 2a(\alpha)$ and $2a(\alpha) \geq \sum_{i = 1}^rk_i(\alpha)$. Consider the composition
\[
 \phi_{a(\alpha), \bold k(\alpha)} \circ \psi_\alpha : \Sec(\mathcal X/B, \alpha) \to Z_{a(\alpha), \bold k(\alpha)}^{\circ} \to U_{\bold k(\alpha)}.
\]
Then 
\[
\phi_{a(\alpha), \bold k(\alpha)} (\psi_\alpha(\Sec(\mathcal X/B, \alpha))) \subset U_{\bold k},
\]
is Zariski open.
Moreover the dimension of $$U_{\bold k(\alpha)} \setminus  \phi_{a(\alpha), \bold k(\alpha)} (\psi_\alpha(\Sec(\mathcal X/B, \alpha))),$$ is at most 
\[
\max\left\{\sum_i k_i (\alpha) - (2a(\alpha) +1-\sum_i k_i (\alpha)), \sum_i k_i (\alpha) -(h(\alpha) - 2a(\alpha) + 1) \right\}.
\]
\end{prop}

Let $C \in \Sec(\mathcal X/B, \alpha)$ and $C' = \rho_F(C)$. Then all members of $\psi_\alpha^{-1}(\psi_\alpha(C))$ are linearly equivalent to each other in $Y_{C'}$.
Thus we can partially compactify $\psi_\alpha : \Sec(\mathcal X/B, \alpha) \to Z_{a(\alpha), \bold k(\alpha)}^{\circ}$ to the relative linear system $$\overline{\psi}_\alpha : \mathcal P_\alpha \to Z_{a(\alpha), \bold k(\alpha)}^{\circ}.$$

\begin{lemm}[{\cite[Theorem 4.10]{DLTT25}}]
\label{lemm:projectivebundleforpsi}
Assume that $h(\alpha) \geq 2a(\alpha)$.
The structure morphism 
\[
\overline{\psi}_\alpha :  \overline{\psi}_\alpha^{-1}(\psi_\alpha( \Sec(\mathcal X/B, \alpha)))  \to\psi_\alpha( \Sec(\mathcal X/B, \alpha)),
\]
is a $\mathbb P^{h(\alpha)-2a(\alpha) + 1}$-bundle over the Zariski open subset $\psi_\alpha( \Sec(\mathcal X/B, \alpha))$. Moreover it is locally trivial in Zariski topology.
\end{lemm}

As a conclusion, we have the following corollary:

\begin{coro}
Assume that $h(\alpha) \geq 2a(\alpha)$ and $2a(\alpha) - \sum_{i = 1}^rk_i(\alpha) \geq 0$. Then $\Sec(\mathcal X/B, \alpha)$ is smooth.
\end{coro}

\section{Rationality of the space of rational curves}

In this section, we establish rationality of $\mathrm{Mor}(\mathbb P^1, S, \alpha)$ for certain nef classes $\alpha \in \mathrm{Nef}_1(S)_{\mathbb Z}$ assuming that $S$ is split, i.e., $\rho(S) = \rho(\overline{S})$. We assume that our ground field is an arbitrary field $k$. 
Let $S$ be a smooth del Pezzo surface defined over $k$. Let $\alpha \in \mathrm{Nef}_1(S)_{\mathbb Z}$ be a nef class on $S$. We fix a blow down $\phi : S \to S_0 \cong \mathbb P^1 \times \mathbb P^1$ such that
\[
2a(\alpha) - \sum_{i = 1}^r k_i(\alpha) \geq 0, \quad 2a'(\alpha) - \sum_{i = 1}^r k_i(\alpha) \geq 0.
\]
The first lemma in this section is the following

\begin{lemm}
\label{lemm:configurationcover_rational}
Assume that we have $2a(\alpha) \geq \sum_{i = 1}^r k_i(\alpha)$. Then the configuration cover 
\[
Z_{a(\alpha), \bold k(\alpha)}^{\circ},
\]
is an irreducible rational variety.
\end{lemm}

\begin{proof}
It follows from Lemma~\ref{lemm:configurationcover_fibration} that the morphism
\[
\phi_{a(\alpha), \bold k(\alpha)} : Z_{a(\alpha), \bold k(\alpha)}^{\circ} \to U_{\bold k(\alpha)},
\]
is dominant and the fibration over the image is a Zariski open subset of a projectivization of a locally free sheaf.
Since $U_{\bold k(\alpha)}$ is a Zariski open subset of the product of projective spaces,
our assertion follows.
\end{proof}

\begin{theo}
\label{theo:rationality}
Let $S$ be a smooth del Pezzo surface of degree $e$ defined over $k$.
When the characteristic $p$ of $k$ is positive, we assume that either $e \geq 2$ or $e = 1$ and $p \geq 11$. Furthermore we may assume that our del Pezzo surface $S$ is not isomorphic to one of surfaces specified in \cite[Theorem 1.2]{BLRT21}.

Let $\alpha \in \mathrm{Nef}_1(S)_{\mathbb Z}$ be a non-zero nef class.
Assume that $h(\alpha) \geq 2a(\alpha)$ and $2a(\alpha) \geq \sum_{i = 1}^r k_i(\alpha)$ hold.
Then the morphism scheme
\[
\mathrm{Mor}(\mathbb P^1, S, \alpha),
\]
is rational.
\end{theo}

\begin{proof}
First it follows from Lemma~\ref{lemm:configurationcover_rational} that $Z_{a(\alpha), \bold k(\alpha)}^{\circ}$ is an irreducible rational variety.
Then by Proposition~\ref{prop:dominant}, the morphism
\[
\Phi_\alpha : \mathrm{Mor}(\mathbb P^1, S, \alpha) \to Z_{a(\alpha), \bold k(\alpha)}^{\circ},
\]
is dominant. 
Then, by Lemma~\ref{lemm:projectivebundleforpsi}, its restriction
\[
\psi_\alpha |_{\psi_\alpha^{-1}(\psi_\alpha(\mathrm{Mor}(\mathbb P^1, S, \alpha)))} : \psi_\alpha^{-1}(\psi_\alpha(\mathrm{Mor}(\mathbb P^1, S, \alpha))) \to \psi_\alpha(\mathrm{Mor}(\mathbb P^1, S, \alpha)),
\]
is a Zariski open subset of the projectivization of a locally free sheaf over $\psi_\alpha(\mathrm{Mor}(\mathbb P^1, S, \alpha))$. Thus our assertion follows.
\end{proof}

\begin{proof}[Proof of Theorem~\ref{theo:intromainIV}]
It follows from Lemma~\ref{lemm:identifyingconicfibration} and \cite[Proposition 3.3]{DLTT25} that one can find a birational morphism $\phi : S \to \mathbb P^1 \times \mathbb P^1$ satisfying the assumptions of Theorem~\ref{theo:rationality}. Thus our assertion follows from Theorem~\ref{theo:rationality}.
\end{proof}

One outstanding question is:
\begin{ques}
Are these moduli spaces still rational when $S$ is not split?
\end{ques}

It would be interesting to explore this question in the view of unirationality of del Pezzo surfaces.

\section{Our set up}
\label{sec:setup}

We briefly recall our set up from \cite[Section 6]{DLTT25}. We assume that our ground field $k$ is a perfect field. Let $S_0$ be $\mathbb P^1 \times \mathbb P^1 = \mathbb P \times \mathbb P'$. Let $\phi : S \to S_0$ be the blow up at $r$ distinct $k$-rational points $(p_i, p_i')$ for $i = 1, \cdots, r$. We assume that $p_i$'s are distinct on $\mathbb P$ as well as $p_i'$'s on $\mathbb P'$.
Let $E_i$ be the exceptional divisor at $(p_i, p_i')$ and we also let $F, F'$ be general fibers of two conic fibrations $S \to S_0 \to \mathbb P^1$.

Let $\alpha$ be a nef class on $S$. Recall that we define the following invariants:
\[
h(\alpha) = -K_S. \alpha, \, a(\alpha) = F.\alpha, \, a'(\alpha) = F'.\alpha, \, k_i(\alpha) = E_i.\alpha.
\]
Throughout the paper we impose the following inequalities:
\[
2a(\alpha) \geq \sum_{i = 1}^r k_i(\alpha), \, 2a'(\alpha) \geq \sum_{i = 1}^r k_i(\alpha).
\]
To ease the notation we denote the moduli space of rational curves $\mathrm{Mor}(\mathbb P^1, S, \alpha)$ by $M_\alpha$.

Let $B = \mathbb P^1$ as usual. Let $f : B \to S$ be a rational curve of the class $\alpha$. This rational curve corresponds to a rational curve $\phi \circ f : B \to S_0$ of the class $(a(\alpha), a'(\alpha))$ such that the multiplicity of this rational curve at $(p_i, p_i')$ is given by $k_i(\alpha)$. Conversely if we have such a rational curve $f' : B \to S_0$, then taking the strict transform, the resulting curve is a rational curve of class $\alpha$.

\subsection{The space of sections}
\label{subsec:spaceofsections}
The above observation motivates the following construction of certain spaces of sections: let $\mathbb P \times \mathbb P' = \mathbb P(V_1) \times \mathbb P(V_2)$ where $V_i$ is a $2$-dimensional vector space over $k$.
Let $V = V_1 \oplus V_2$, $\mathcal V_1 = \mathcal O(a) \otimes V_1, \mathcal V_2 = \mathcal O(a') \otimes V_2, \mathcal V = \mathcal V_1 \oplus \mathcal V_2$ be vector bundles over $B = \mathbb P^1$ where $a, a'$ be non-negative integers. Let $\ell_{i, j} \subset V_j$ be the $1$-dimensional subspaces corresponding to $p_i, p_i'$ respectively.
Let $\widetilde{M}_{\alpha} = \widetilde{M}_{a, a', \bold k}$ be the space of sections
\[
(s, t) \in H^0(B, \mathcal V_1) \oplus H^0(B, \mathcal V_2),
\]
such that
\begin{itemize}
\item both $s, t$ are nowhere vanishing, and;
\item the length of $(s, t)^{-1}(\ell_{i, 1} \oplus \ell_{i, 2})$ is $k_i$ for $i = 1, \cdots, r$ where $k_i$ is a non-negative integer.
\end{itemize}
Then this scheme admits a natural morphism
\[
\widetilde{M}_\alpha \to M_\alpha,
\]
which realizes $\widetilde{M}_\alpha$ as a $\mathbb G_m^2$-torsor over $M_\alpha$.

Let $U_{\bold k}$ as before. For any $(w_i) \in U_{\bold k}$, we define
\[
E_w \subset H^0(B, \mathcal V),
\]
as the space of sections $(s, t)$ such that $$(s, t)(w_i) \subset \mathcal O(a)\otimes\ell_{i, 1} \oplus \mathcal O(a')\otimes \ell_{i, 2}.$$ This is bundled as $E \to U_{\bold k}$ and we realize $\widetilde{M}_\alpha \to M_\alpha \to U_{\bold k}$ as a Zariski open subset of $E \to U_{\bold k}$.

\subsection{Stratifications by poschemes}

Point counting on $E$ is not so difficult as this is generically a vector bundle over $U_{\bold k}$, and its point counting is reduced to point counting on $U_{\bold k}$. Thus our goal is to understand the complement
\[
E \setminus \widetilde{M}_{a, a', \bold k}.
\]
To understand this object, Das--Tosteson introduced certain simplicial resolutions in \cite{DT24}. We recall their constructions in \cite{DT24} and \cite{DLTT25}, and we freely use the notation from these papers. Let $Q$ be a finite \'etale poscheme (a scheme with a poset structure as \cite[Definition 3.1.1]{DH24}) consisting of subspaces of $V$:
\[
Q =\{ V, V_1, V_2, \ell_{i, 1}\oplus \ell_{i, 2}, \ell_{i, 1}, \ell_{i, 2}, \, (i = 1, \cdots, r), 0\},
\]
where the poset structure is given by inclusion. 
Note that this poscheme admits meets, i.e., any intersection of elements in $Q$ is contained in $Q$.
Let $\mathrm{Hilb}(B)$ be the Hilbert scheme parametrizing all subschemes of $B$ and we consider this scheme as a poscheme under inclusion. Then we construct the scheme $\mathrm{Hilb}(B)^Q$ parametrizing homomorphisms
\[
x : Q \to \mathrm{Hilb}(B),
\]
such that $x^{-1}(B) = V$.
Then we embed $w \in U_{\bold k}$ into $\mathrm{Hilb}(B)^Q$ by letting
\[
x_q = 
\begin{cases}
B & \text{ if $q = V$}\\
w_i & \text{ if $q = \ell_{i, 1}\oplus \ell_{i, 2}$}\\
\emptyset & \text{ otherwise.}
\end{cases}
\]
Then we consider the following poscheme over $U_{\bold k}$:
\[
(U_{\bold k} < \mathrm{Hilb}(B)^Q) = \{(w < x) \, | \, w \in U_{\bold k}, x \in \mathrm{Hilb}(B)^Q \text{ such that $w < x$. } \}.
\]
For each $(w < x) \in (U_{\bold k} < \mathrm{Hilb}(B)^Q)$, we define the space $Z_{w < x} \subset E_w$ as the space of $(s, t) \in E_w$ such that $(s, t)(x_q) \subset \mathcal K_q \subset \mathcal V$ for any $q \in Q$ where $\mathcal K_q$ is a subbundle of $\mathcal V$ corresponding to $q \in Q$.
Then these spaces are bundled as a stratification by linear subspaces:
\[
Z \subset (U_{\bold k} < \mathrm{Hilb}(B)^Q) \times_{U_{\bold k}} E,
\]
in the sense of \cite[Definition 2.16]{DLTT25}.
Then the complement is the image of $Z$, i.e., 
\[
E \setminus \widetilde{M}_{a, a', \bold k} = \mathrm{im}(Z \to E).
\]
In this way, $Z$ provides a simplicial resolution of $E \setminus \widetilde{M}_{a, a', \bold k}$ and the method of the bar complexes developed in \cite{DT24, DLTT25} is kicked in.

\subsection{Bar complexes}

An element $x \in \mathrm{Hilb}(B)^Q(\overline{k})$ is called saturated in the sense of \cite[Definition 5.2]{DLTT25} if $x$ preserves meets, i.e., $x_{p\wedge q} = x_p \cap x_q$ for any $p, q \in Q$.
Such elements form a Zariski open subscheme $$Q^{JB} \subset \mathrm{Hilb}(B)^Q.$$
See \cite[Definition 5.2]{DLTT25} for more details.

As \cite[Definition 5.16]{DLTT25}, there is a stratification $$\mathrm{Hilb}(B)^Q = \sqcup_{T'} \mathcal N_{T'},$$ into locally closed subsets parametrized by combinatorial types $T'$ of $\mathrm{Hilb}(B)^Q$ as in \cite[Section 5.4]{DLTT25}. Then our $U_{\bold k} \subset \mathrm{Hilb}(B)^Q$ is a locally closed subset which is a finite union of $\mathcal N_{T'}$ of saturated combinatorial types.

Let $P \subset (U_{\bold k} < \mathrm{Hilb}(B)^Q)$ be a closed and downward closed union of finitely many locally closed set $\mathcal Z_T$ for finitely many saturated types $T$ of $(U_{\bold k} < \mathrm{Hilb}(B)^Q)$ where $\mathcal Z_T$ is the topological, downward, and saturation closure of $\mathcal N_T$. We assume that the projection $P \to U_{\bold k}$ is proper.

Let $U\subset U_{\bold k}$ be a Zariski open subset and we denote the base changes of $(U_{\bold k} < \mathrm{Hilb}(B)^Q)$ and $P$ via $U \to U_{\bold k}$ by $(U < \mathrm{Hilb}(B)^Q)$ and $P_U$ respectively.
Recall that we have a fibration $E \to U_{\bold k}$ and we suppose that for its restriction $E|_U \to U$ every fiber has constant rank. We also have a stratification $Z \subset (U_{\bold k} < \mathrm{Hilb}(B)^Q) \times_{U_{\bold k}} E$ and we denote its restriction to $P_U$ by $Z_{P_U}$.
Then one can define the bar complex $B(P_U, Z_{P_U})$ as in \cite[After Definition 2.17]{DLTT25} which is a certain simplicial scheme. This assigns for $[n] = \{0, \cdots, n\}$, the fiber product:
\[
B(P_U, Z_{P_U})([n]) = \{ w < x_0 \leq \cdots \leq x_n, z \in Z_{w < x_n} \, | \, (w < x_i) \in P_U\},
\]
with an augmentation morphism $B(P_U, Z_{P_U}) \to E|_U$ mapping 
\[
(w < x_0 \leq \cdots \leq x_n, z)\mapsto (w, z).
\]

For any saturated combinatorial type $T$ of $(U_{\bold k} < \mathrm{Hilb}(B)^Q)$, one can define the functions
\[
\gamma(T), \, \mathrm{rank}(T), \, |\mathrm{Supp}(T)|,
\]
which are defined in \cite[Example 5.14]{DLTT25}, \cite[Example 5.12]{DLTT25}, and \cite[Example 5.13]{DLTT25} respectively. The function $\gamma(T)$ measures the expected codimension of $Z_{w < x}$ in $E_w$. The function $\mathrm{rank}(T)$ measures the rank of $x$ in the poset $[w, x]$ in $Q^{JB}(\overline{k})$. The function $|\mathrm{Supp}(T)|$ is the number of points in the support of $x$. Then for any saturated combinatorial type $T$ of $(U_{\bold k} < \mathrm{Hilb}(B)^Q)$, we define the function:
\[
\kappa(T) = 2\gamma(T) - \mathrm{rank}(T) - 2|\mathrm{Supp}(T)|.
\]
One may think of $\kappa(T)$ as the expected homological codimension of the contribution of $\mathcal Z_T$. The following theorem has been proved over $\mathbb C$ in \cite[Theorem 5.9]{DT24} and extended to arbitrary perfect fields in \cite[Theorem 7.5]{DLTT25}:
\begin{theo}[{\cite[Theorem 5.9]{DT24} and \cite[Theorem 7.5]{DLTT25}}]
\label{theo:barcomplex}
Let $I \in \mathbb Q_{\geq 0}$. Assume that
\begin{itemize}
\item $P \to U_{\bold k}$ is proper and $P$ is closed and downward closed in $(U_{\bold k} < \mathrm{Hilb}(B)^Q)$;
\item for any $(w < x) \in P_U(\overline{k})$ with $x$ being saturated and every $y \in Q^{JB}(\overline{k})$ such that $x < y$ and $(x, y)\cap Q^{JB}(\overline{k}) = \emptyset$, the space $Z_{w < y} \subset E_w$ has expected codimension $\gamma(w < y)$, and;
\item $P$ contains every $\mathcal Z_T$ such that $\kappa(T) \leq I$.
\end{itemize}
Then for all $i > 2\dim E|_U - \lfloor I\rfloor -2$, the augmentation morphism
\[
B(P_U, Z_{P_U}) \to \mathrm{im}(Z_{P_U} \to E|_U),
\]
induces an isomorphism
\[
H^i_{\text{\'et}, c}(\mathrm{im}(Z_{P_U} \to E|_U)_{\overline{k}}, \mathbb Q_\ell) \cong H^i_{\text{\'et}, c}(B(P_U, Z_{P_U})_{\overline{k}}, \mathbb Q_\ell),
\]
as Galois representations of $\ell$-adic cohomologies.
\end{theo}

Thus to understand the cohomology of 
\[
E \setminus \widetilde{M}_{a, a', \bold k} = \mathrm{im}(Z \to E),
\]
it suffices to analyze the cohomology of the bar complex $B(P_U, Z_{P_U})$.

\section{Main results}
\label{sec:Manin}

In this section, we prove the main result of this paper which follows from the discussion of \cite[Section 8]{DLTT25}. We closely follow its exposition. We assume that our ground field $k$ is a finite field $\mathbb F_q$.
We work on the set up introduced in Section~\ref{sec:setup} so that $S$ is a blow up of $S_0 = \mathbb P^1 \times \mathbb P^1$ at $r$ $k$-rational points. Let $\alpha$ be a nef class on $S$ with invariants:
\[
a(\alpha) = F.\alpha, \, a'(\alpha) = F'.\alpha, \, E_i.\alpha = k_i(\alpha),
\]
where $F, F', E_i$ are divisors defined in Section~\ref{sec:setup}.

\subsection{Upper bounds}

The following lemma is essentially \cite[Lemma 8.1]{DLTT25}, and it follows from our geometric results established in Section~\ref{sec:geometryofmoduli}:
\begin{lemm}
\label{lemm:upperbounds}
Assume that $S$ is split and the following inequalities hold:
\[
2a(\alpha) \geq \sum_{i = 1}^r k_i(\alpha), \, 2a'(\alpha) \geq \sum_{i = 1}^r k_i(\alpha).
\]
Then we have
\[
\#M_\alpha(k) \leq \frac{q^{h(\alpha) + 2}}{(1-q^{-1})^{r+2}}.
\]
\end{lemm}

\begin{proof}
We include a proof for completeness.
First we consider the morphism
\[
\psi_\alpha : M_\alpha \to Z^\circ_{a(\alpha), \bold k(\alpha)},
\]
constructed in Section~\ref{subsec:fibrationstr}.
For each $k$-rational point $x$ on $\psi_\alpha(M_\alpha)$, it follows from Lemma~\ref{lemm:projectivebundleforpsi} that the fiber $\psi_\alpha^{-1}(x)$ is a Zariski open subset of the projective space of dimension $h(\alpha) - 2a(\alpha) + 1$. This implies that
\[
\#M_\alpha(k)  \leq \#Z^\circ_{a(\alpha), \bold k(\alpha)}(k)\frac{q^{h(\alpha) - 2a(\alpha) + 1}}{(1-q^{-1})}.
\]
Next we consider the morphism
\[
\phi_{a(\alpha), \bold k(\alpha)} : Z^\circ_{a(\alpha), \bold k(\alpha)} \to U_{\bold k(\alpha)},
\]
constructed in Section~\ref{subsec:furtherZ}. It follows from Lemma~\ref{lemm:configurationcover_fibration} that for any $k$-rational point $x$ on $\phi_{a(\alpha), \bold k(\alpha)} (Z^\circ_{a(\alpha), \bold k(\alpha)} )$, the fiber $\phi_{a(\alpha), \bold k(\alpha)}^{-1}(x)$ is a Zariski open subset of the projective space of dimension $2a(\alpha) +1 -\sum_{i = 1}^r k_i(\alpha)$. This implies
\[
\#Z^\circ_{a(\alpha), \bold k(\alpha)}(k) \leq \#U_{\bold k(\alpha)} (k) \frac{q^{2a(\alpha) + 1 - \sum_i k_i(\alpha)}}{(1-q^{-1})}.
\]
Finally since $U_{\bold k(\alpha)}$ is a Zariski open subset of $\prod_{i = 1}^r \mathbb P^{k_i}$, we have
\[
\#U_{\bold k(\alpha)} (k) \leq \frac{q^{\sum_{i = 1}^rk_i(\alpha)}}{(1-q^{-1})^r}.
\]
Combining everything together our assertion follows.

\end{proof}

Using the above lemma, one can obtain upper bounds of correct magnitude for split del Pezzo surfaces of degree $\geq 4$. This has been proved in \cite[Theorem 8.3]{DLTT25} for split quartic del Pezzo surfaces, and the same proof applies to split quintic del Pezzo surfaces. We record it here for completeness:

\begin{theo}
\label{theo:upperbound}
Let $k = \mathbb F_q$, $B = \mathbb P^1$, and $S$ be a split smooth del Pezzo surface of degree $5$ defined over $k$.
We define the counting function:
\[
N(B, S, -K_S, d) = \sum_{\alpha \in \mathrm{Nef}_1(S)_{\mathbb Z}, \, -K_S.\alpha \leq d} \#M_\alpha(k).
\]
Then we have
\[
\limsup_{d \to \infty} \frac{N(B, S, -K_S, d)}{q^dd^4} \leq \frac{\alpha(-K_S, \mathrm{Nef}_1(S))q^2}{(1-q^{-1})^6},
\]
where $\alpha(-K_S, \mathrm{Nef}_1(S))$ is the alpha constant of the cone $\mathrm{Nef}_1(S)$. See \cite[Definition 8.2]{DLTT25} for its definition.
\end{theo}

\begin{proof}
This follows from Lemma~\ref{lemm:upperbounds}. See the proof of \cite[Theorem 8.3]{DLTT25} for more details.
\end{proof}

\subsection{Counting curves through the bar complexes}
\label{subsec:countingthroughbar}

Our goal of this section is to describe $\#\widetilde{M}_{a, a', \bold k}(k)$ in a combinatorial way in terms of the poset $(U_{\bold k} \leq Q^{JB})(k)$ up to error terms, and our discussion is heavily based on \cite[Section 8.3]{DLTT25}. We closely follow its exposition.
Suppose that we have positive integers $h, a, a', k_i (i = 1, \cdots, r)$ such that
\[
h = 2a + 2a' - \sum_{i = 1}^rk_i, \, 2a - \sum_{i = 1}^rk_i >0, \, 2a' - \sum_{i = 1}^rk_i >0.
\]
As in \cite[Section 8.3]{DLTT25}, we set 
\[
I = \frac{1}{8}\min\left\{2a - \sum_{i = 1}^rk_i, 2a' - \sum_{i = 1}^rk_i \right\} -\frac{1}{2}.
\]
For any $q \in Q$, let $m_q : \mathrm{Hilb}(B)^Q(\overline{k}) \to \mathbb N$ be the function defined in \cite[Example 5.11]{DLTT25}. For any $(w < x) \in (U_{\bold k} < \mathrm{Hilb}(B)^Q)(\overline{k})$, we define
\begin{align*}
&E(w < x) := \gamma(w < x)  \\
&= 4m_0(x) + 3\sum_{i, j}m_{\ell_{i, j}}(x) + 2\sum_j m_{V_j}(x) + 2\sum_{i}m_{\ell_{i, 1}\oplus \ell_{i, 2}}(x) - 2\sum_i k_i.
\end{align*}
Let $P \subset (U_{\bold k} < \mathrm{Hilb}(B)^Q)$ be the closed subposcheme consisting of $(w < x)$ such that $E(w < x) \leq 2I$. Then it follows from \cite[Lemma 8.4]{DLTT25} that $P$ is proper over $U_{\bold k}$. Moreover since coefficients of $E$ are decreasing, one can conclude that $P$ is downward closed and it consists of finitely many $\mathcal N_T$.

Now we consider the morphism
\[
M_{a, a', \bold k} \to U_{\bold k},
\]
constructed in Section~\ref{subsec:fibrationstr}.
Following the same construction of \cite[Section 8.3, p. 45--46]{DLTT25}, 
we obtain a proper closed subset $$F = F_1\cup F_2\cup F_3 \subsetneq U_{\bold k}.$$ 
Note that the only difference from \cite[Section 8.3, p. 45--46]{DLTT25} is that the number of singular fibers is given by $r = r(S)$, but not $4$. This slightly changes the numerology, but it is almost identical. In particular, letting $U = U_{\bold k} \setminus F$, the assumptions of Theorem~\ref{theo:barcomplex} are satisfied.
Then the fact that the codimension of $F$ in $U_{\bold k}$ is greater than or equal to $4I$ follows from  Proposition~\ref{prop:codimensionI} instead of \cite[Proposition 4.13]{DLTT25}.

Let $\widetilde{M}_{a, a', \bold k} \to M_{a, a', \bold k}$ be the $\mathbb G_m^2$-torsor considered in Section~\ref{subsec:spaceofsections}. Our goal is to understand
\[
\# M_{a, a', \bold k}(k) = (q-1)^{-2}\# \widetilde{M}_{a, a', \bold k}(k).
\] 
It follows from the Grothendieck--Lefschetz trace formula that
\[
\# \widetilde{M}_{a, a', \bold k}(k) = \sum_i (-1)^i \mathrm{Tr}(\mathrm{Frob}\curvearrowright H^i_{\text{\'et}, c}((M_{a, a', \bold k})_{\overline{k}}, \mathbb Q_\ell)).
\]
Then by Theorem of Sawin--Shusterman (\cite[Theorem B.1]{DLTT25}), there exists a uniform constant $C_1$ such that 
\[
\sum_{i < 4a + 4a' -2\sum_i k_i -I -10}(-1)^i \mathrm{Tr}(\mathrm{Frob}\curvearrowright H^i_{\text{\'et}, c}((M_{a, a', \bold k})_{\overline{k}}, \mathbb Q_\ell)) = O(q^{2a + 2a' - \sum_ik_i - I/2+5}C_1^h).
\]

Let $U = U_{\bold k} \setminus F$. Since the codimension of $F$ in $U_{\bold k}$ is greater than or equal to $4I$ and $\widetilde{M}_{a, a', \bold k} \to U_{\bold k}$ is flat, we have
\begin{align*}
&\sum_{i \geq 4a + 4a' -2\sum_i k_i -I -10}(-1)^i \mathrm{Tr}(\mathrm{Frob}\curvearrowright H^i_{\text{\'et}, c}((M_{a, a', \bold k})_{\overline{k}}, \mathbb Q_\ell))\\
&=\sum_{i \geq 4a + 4a' -2\sum_i k_i -I -10}(-1)^i \mathrm{Tr}(\mathrm{Frob}\curvearrowright H^i_{\text{\'et}, c}((M_{a, a', \bold k}|_U)_{\overline{k}}, \mathbb Q_\ell)).
\end{align*}
Since we have $E|_U = M_{a, a', \bold k}|_U \sqcup \mathrm{im}(Z_{P_U} \to E|_U)$, as discussed in \cite[Section 8.3]{DLTT25}, we have
\begin{align*}
&\sum_{i \geq 4a + 4a' -2\sum_i k_i -I -10}(-1)^i \mathrm{Tr}(\mathrm{Frob}\curvearrowright H^i_{\text{\'et}, c}((M_{a, a', \bold k}|_U)_{\overline{k}}, \mathbb Q_\ell))\\
&= \sum_{i \geq 4a + 4a' -2\sum_i k_i -I -10}(-1)^i \mathrm{Tr}(\mathrm{Frob}\curvearrowright H^i_{\text{\'et}, c}((E|_U)_{\overline{k}}, \mathbb Q_\ell))\\
&\qquad - \sum_{i \geq 4a + 4a' -2\sum_i k_i -I -10}(-1)^i \mathrm{Tr}(\mathrm{Frob}\curvearrowright H^i_{\text{\'et}, c}(\mathrm{im}(Z_{P_U} \to E|_U)_{\overline{k}}, \mathbb Q_\ell))\\
&\qquad \qquad + O(q^{2a + 2a' - \sum_ik_i - I/2+6}C_1^h).
\end{align*}
We are interested in
\[
\sum_{i \geq 4a + 4a' -2\sum_i k_i -I -10}(-1)^i \mathrm{Tr}(\mathrm{Frob}\curvearrowright H^i_{\text{\'et}, c}(\mathrm{im}(Z_{P_U} \to E|_U)_{\overline{k}}, \mathbb Q_\ell)).
\]
It follows from Theorem~\ref{theo:barcomplex} that the above quantity is equal to
\[
\sum_{i \geq 4a + 4a' -2\sum_i k_i -I -10}(-1)^i \mathrm{Tr}(\mathrm{Frob}\curvearrowright H^i_{\text{\'et}, c}(B(P_U, Z_{P_U})_{\overline{k}}, \mathbb Q_\ell)).
\]
Applying the discussion of \cite[Section 8.3 p. 48--50]{DLTT25} we conclude
\begin{align*}
&  \sum_{i \geq 4a + 4a' -2\sum_i k_i -I -10}(-1)^i \mathrm{Tr}(\mathrm{Frob}\curvearrowright H^i_{\text{\'et}, c}(\mathrm{im}(Z_{P_U} \to E|_U)_{\overline{k}}, \mathbb Q_\ell))\\
& = -q^{2a +2a' + 4}\sum_{(w < x) \in (U_{\bold k} < Q^{JB})(k)} \mu_k(w, x)q^{-\gamma(x)}
+O(q^{2a + 2a' - \sum_i k_i -I/4 + 6} \max\{C_1, C_2\}^h),
\end{align*}
where $C_2$ is a constant coming from \cite[Theorem 7.9]{DLTT25} with $\eta = 1/2$ and $$\mu_k : (Q^{JB}(k) \leq Q^{JB}(k))\to \mathbb Z,$$ is the M\"obius function for the poset $Q^{JB}(k)$. (See \cite[Definition 7.6]{DLTT25} for the definition.)
Similarly one can conclude
\begin{align*}
&  \sum_{i \geq 4a + 4a' -2\sum_i k_i -I -10}(-1)^i \mathrm{Tr}(\mathrm{Frob}\curvearrowright H^i_{\text{\'et}, c}((E|_U)_{\overline{k}}, \mathbb Q_\ell))\\
& = q^{2a +2a' + 4}\sum_{w  \in U_{\bold k}(k)} q^{-\gamma(w)}
+O(q^{2a + 2a' - \sum_i k_i -I/4 + 6} \max\{C_1, C_2\}^h).
\end{align*}
Combining everything together we conclude
\begin{theo}
\label{theo:maintheoremI}
Assume that $q > C_3$. Then we have
\begin{align*}
&\# M_{a, a', \bold k}(k) = \\
&q^{2a +2a' + 4}\sum_{(w \leq x) \in (U_{\bold k} \leq Q^{JB})(k)} \mu_k(w, x)q^{-\gamma(x)}
+O(q^{2a + 2a' - \sum_i k_i -I/4 + 6} \max\{C_1, C_2\}^h).
\end{align*}
\end{theo}

\subsection{Virtual curve counting}
We follow the exposition of \cite[Section 8.4]{DLTT25}. The difference in this section from \cite[Section 8.4]{DLTT25} is that certain formulas are generalized from $r = 4$ to arbitrary $r$. 
For this reason, we record the modified formulas, but do not include detailed proofs.

Our goal here is to understand the asymptotic behavior of 
\[
q^{2a +2a' + 4}\sum_{(w \leq x) \in (U_{\bold k} \leq Q^{JB})(k)} \mu_k(w, x)q^{-\gamma(x)},
\]
as $a, a', k_i \to \infty$ under the conditions 
\[
2a - \sum_{i = 1}^r k_i >0, \, 2a - \sum_{i = 1}^r k_i >0.
\]
As \cite[Section 8.4]{DLTT25}, we consider the following virtual height zeta function:
\begin{align*}
&\mathsf Z(\bold t)\\
&=\sum_{i = 1, \cdots, r, \, k_i \in \mathbb N} q^{\sum_{i = 1}^rk_i} \left( \sum_{(w \leq x) \in (U_{\bold k} \leq Q^{JB})(k)} \mu_k(w, x)q^{-\gamma(x)}\right) \prod_{i = 1}^r t_i^{k_i}.
\end{align*}
Since the M\"obius function is multiplicative as \cite[Lemma 7.7(1)]{DLTT25}, the above zeta function becomes the following convergent Euler product:
\begin{align*}
&\prod_{c \in |B|} \Bigg(\sum_{x \in \mathrm{ch}(Q)} \mu_k(\hat{1}, x)q^{-\gamma(x)|c|} \\
&\qquad \qquad + \sum_{i = 1}^r\sum_{d_i = 1}^\infty(qt_i)^{|c|d_i}\sum_{x \in \mathrm{ch}(Q), d_i[\ell_{i, 1}\oplus\ell_{i, 2}] \leq x}
\mu_k(d_i[\ell_{i, 1}\oplus\ell_{i, 2}], x)q^{-\gamma(x)|c|} \Bigg),
\end{align*}
where $|B|$ is the set of closed points on $B$, $|c|$ is $[k(c):k]$, $\mathrm{ch}(Q)$ is the set of chains as \cite[Definition 5.3]{DLTT25}, and $d[\ell_{i, 1}\oplus\ell_{i, 2}]$ is a certain chain defined in \cite[Example 5.5]{DLTT25}. Then it follows from \cite[Lemma 7.7(2) and Example 5.8]{DLTT25} that this Euler product becomes
\begin{align*}
&\prod_{c \in |B|} \Bigg(1 -(r+2)q^{-2|c|} + 2rq^{-3|c|} -(r-1)q^{-4|c|} \\
&\qquad  + \sum_{i = 1}^r\sum_{d_i = 1}^\infty(qt_i)^{|c|d_i}(q^{-2d_i|c|} -2q^{-(2d_i+1)|c|} +2q^{-(2d_i+3)|c|} - q^{-(2d_i+4)|c|})\Bigg).
\end{align*}
This expression shows that this Euler product absolutely converges when $|t_i|\leq q^{-\delta}$ for $i = 1, \cdots, r$ for any fixed $\delta > 0$.
\begin{prop}[{\cite[Proposition 8.7]{DLTT25}}]
We have the following Euler product:
\begin{align*}
&\prod_{i = 1}^rZ_{\mathbb P^1}(q^{-1}t_i)\mathsf Z(\bold t)  \\
&=\prod_{c \in |B|} \Bigg(1 + O(q^{-2|c|})(1 + \sum_{i = 1}^r t_i^{|c|}) + \sum_{l_i \in \{0, 1\}, \sum_i l_i > 1}O(q^{-(\sum_il_i)|c|})\prod_{i = 1}^r t_i^{l_i|c|} \Bigg),
\end{align*}
where $Z_{\mathbb P^1}(t)$ is the Hasse-Weil zeta function for $\mathbb P^1$. In particular there exists a uniform $\delta > 0$ such that this Euler product absolutely converges when $|t_i|\leq q^\delta$ for $i = 1, \cdots, r$. Moreover, we have
\[
\lim_{t_i \to 1}\prod_{i = 1}^r (1-t_i)\mathsf Z(\bold t) = (1-q^{-1})^{-r}\prod_{c \in |B|} (1-q^{-|c|})^{r+2}(1 +(r+2)q^{-|c|} + q^{-2|c|} ).
\]
\end{prop}

\begin{proof}
One may argue as \cite[Proposition 8.7]{DLTT25}.
\end{proof}

Everything put together and combining them with Abel's summation techniques, we conclude
\begin{prop}[{\cite[Proposition 8.8]{DLTT25}}]
\label{prop:main}
There exists $\eta > 0$ which does not depend on $q, \bold k$ such that
\begin{align*}
&q^{\sum_{i = 1}^rk_i} \left( \sum_{(w \leq x) \in (U_{\bold k} \leq Q^{JB})(k)} \mu_k(w, x)q^{-\gamma(x)}\right)\\
&=  (1-q^{-1})^{-r}\prod_{c \in |B|} (1-q^{-|c|})^{r+2}(1 +(r+2)q^{-|c|} + q^{-2|c|} ) + O(q^{-\eta\min\{k_i\}}).
\end{align*}
\end{prop}

\begin{proof}
One may argue as \cite[Proposition 8.8]{DLTT25}.
\end{proof}

\subsection{The main result}
Recall that we are in the set up described in Section~\ref{sec:setup}.
Let $\alpha$ be an ample class on $S$ such that 
\[
2a(\alpha) - \sum_{i = 1}^r k_i(\alpha) > 0, \, 2a'(\alpha) - \sum_{i = 1}^r k_i(\alpha) > 0.
\]
Let $\epsilon > 0$ be a sufficiently small number such that
\[
\min\left\{2a(\alpha) - \sum_{i = 1}^r k_i(\alpha), 2a'(\alpha) - \sum_{i = 1}^r k_i(\alpha) \right\} \geq -\epsilon K_S.\alpha.
\]
Here is the main theorem of this section:
\begin{theo}
\label{theo:mainI}
Assume that $q^{\frac{\epsilon}{32}} > \max\{C_1, C_2\}$ as well as $q > C_3$ where $C_1, C_2, C_3$ are constants coming from Section~\ref{subsec:countingthroughbar}. Then we have
\[
\lim_{m \to \infty}\frac{\#M_{m\alpha}(k)}{q^{-mK_S.\alpha}} = q^2 (1-q^{-1})^{-r-2}\prod_{c \in |B|} (1-q^{-|c|})^{r+2}(1 +(r+2)q^{-|c|} + q^{-2|c|} ).
\]
\end{theo}
Note that this limit coincides with Tamagawa number $\tau_{-K_S}(S)$ introduced in \cite{Peyre}.

\begin{proof}
This follows from Theorem~\ref{theo:maintheoremI} and Proposition~\ref{prop:main}.
Note that we have
\[
-K_S.\alpha = h(\alpha) = 2a(\alpha) + 2a'(\alpha) - \sum_{i = 1}^r k_i(\alpha),
\]
and
\[
I = \frac{1}{8}\min\left\{2a(\alpha) - \sum_{i = 1}^rk_i(\alpha), 2a'(\alpha) - \sum_{i = 1}^rk_i(\alpha) \right\} -\frac{1}{2}.
\]
Thus Theorem~\ref{theo:maintheoremI} and Proposition~\ref{prop:main} shows that
\[
\#M_{\alpha}(k)/q^{h(\alpha)} = \tau_{-K_S}(S) + O(\max\{q^{-\epsilon h(\alpha)/32}C^{h(\alpha)}, q^{-\eta\min_i\{E_i.\alpha\}}\}),
\]
where $C = \max\{C_1, C_2\}$.
Assuming $q^{\epsilon/32} > C$, our assertion follows.
\end{proof}

\section{Applications}

Here we collect concrete applications of results in Section~\ref{sec:Manin} and prove main theorems in the introduction.

\subsection{Split quintic del Pezzo surfaces}
\label{subsec:quintic}

Here we apply our results to split quintic del Pezzo surfaces. 
The first easy consequence of Theorem~\ref{theo:mainI} is the following:
\begin{theo}
\label{theo:mainII}
There exist uniform constants $C \geq 1$ and $C_3$ with the following properties:
let $S$ be a smooth split del Pezzo surface of degree $5$ defined over $\mathbb F_q$.
Let $\alpha$ be an ample class on $S$.
By Lemma~\ref{lemm:identifyingconicfibration}, one can find a birational morphism $\phi : S \to S_0 = \mathbb P^1 \times \mathbb P^1$ such that we have
\[
\mathfrak I(\alpha) = \min\left\{ 2F.\alpha - \sum_{i = 1}^3 E_i.\alpha, 2F.'\alpha - \sum_{i = 1}^3 E_i.\alpha \right\} >0,
\]
where $E_1, E_2, E_3$ are exceptional divisors for $\phi$ and $F, F'$ are general fibers of two conic fibrations coming from $S_0$.
Let $\epsilon(\alpha) = \mathfrak I(\alpha)/(-32K_S.\alpha)$.
Suppose that we have $q^{\epsilon(\alpha)} > C$ and $q > C_3$. Then we have
\[
\lim_{m \to \infty}\frac{\#\mathrm{Mor}(\mathbb P^1, S, m\alpha)(k)}{q^{-mK_S.\alpha}} = q^2 (1-q^{-1})^{-5}\prod_{c \in |B|} (1-q^{-|c|})^{5}(1 +5q^{-|c|} + q^{-2|c|} ).
\]
\end{theo}

Next we prove our main theorems in the introduction:

\begin{proof}[Proof of Theorem~\ref{theo:intromainI}]
Let $\mathfrak T$ be the set of birational morphisms $\phi : S \to S_0 = \mathbb P^1 \times \mathbb P^1$. For a nef class $\alpha \in \mathrm{Nef}_1(S)$ and a birational morphism $\phi: S \to S_0$, we define $\mathfrak I_\phi(\alpha)$ as in Theorem~\ref{theo:mainII}. We also define
\[
\mathfrak J_\phi (\alpha) = \min\{ E_1.\alpha, E_2.\alpha, E_3. \alpha\}.
\]
Then one can define a piecewise linear function $\ell(\alpha)$ as
\[
\ell(\alpha) = \max_{\phi \in \mathfrak T}\{ \min\{ \mathfrak I_{\phi}(\alpha) /32, \mathfrak J_\phi(\alpha)\} \}.
\]
Then Theorem~\ref{theo:intromainI} follows from Theorem~\ref{theo:mainII}.
\end{proof}
\begin{proof}[Proof of Theorem~\ref{theo:intromainII}]
It follows from Theorem~\ref{theo:maintheoremI} and Proposition~\ref{prop:main} that for any $\alpha \in \mathrm{Nef}_1(S)_{\epsilon, \mathbb Z}$ we have
\[
\#M_\alpha(k)/q^{-K_S.\alpha} = \tau_{-K_S}(S) + O(\max\{q^{-\epsilon h(\alpha)}C^{h(\alpha)}, q^{-\eta\epsilon h(\alpha)}\}).
\]
Combining this estimate with the argument of \cite[The proof of Theorem 1.1]{DLTT25} using the Ehrhart polynomial of $\mathrm{Nef}_1(S)_{\epsilon}$ will prove the statement under the assumption that $q^\epsilon > C$.
\end{proof}

\begin{rema}
Let us estimate $C_1, C_2, C_3$. Since any smooth degree $5$ del Pezzo surface is defined by five quadrics in $\mathbb P^5$, the proof of \cite[Theorem B.1]{DLTT25} will show that $C_1 = 2^{48}$ works. The constants $C_2, C_3$ come from \cite[Theorem 7.9]{DLTT25}. In its notation, $c_1 = 6$ works, so for $C_2$ one may define $C_2 = 240$ with $\eta = 1/2$. Then $C_3 = (240)^4$ works.
Thus one might estimate $C = 2^{48}$.
\end{rema}

\begin{rema}
It would be nice to establish the asymptotic formula for $N(\mathbb P^1, S, -K_S, \mathrm{Nef}_1(S), d)$. A reason preventing this is that for some of non-ample nef classs $\alpha$, our stability range $\mathfrak I_\phi(\alpha)$ is equal to $0$. This is not enough to beat the exponential error term. We plan to come back to this issue in our future work.
\end{rema}

\subsection{Lower bounds for low degree del Pezzo surfaces}

Finally we prove Theorem~\ref{theo:intromainIII}:
\begin{proof}[Proof of Theorem~\ref{theo:intromainIII}]
We fix a birational morphism $\phi : S \to S_0 = \mathbb P^1 \times \mathbb P^1$ with divisor classes $F, F', E_1, \cdots, E_r$ where $r = 8-e$. Let $\epsilon > 0$ be a small rational number.
We set
\[
\mathcal C = \{ \alpha \in \mathrm{Nef}_1(S)\, | \, \min\{ \mathfrak I(\alpha), \mathfrak J(\alpha)\}/(-K_S.\alpha) \geq \epsilon\},
\]
where
\[
\mathfrak I(\alpha) = \min\left\{2F.\alpha - \sum_{i = 1}^rE_i.\alpha, 2F'.\alpha - \sum_{i = 1}^rE_i.\alpha\right\}
\]
\[
\mathfrak J(\alpha)=  \min\{E_1.\alpha, \cdots, E_r.\alpha\}.
\]
Then a similar proof to Theorem~\ref{theo:intromainII} will prove this theorem.
\end{proof}

\bibliographystyle{alpha}
\bibliography{stability}

\end{document}